\documentstyle[12pt]{article}
\newcommand{\qed}{\hfill\rule{4pt}{8pt}\par\vspace{\baselineskip}}

\newtheorem{de}{Definition}[section]
\newtheorem{lm}[de]{Lemma}

\newtheorem{pr}[de]{Proposition}
\newtheorem{co}[de]{Corollary}
\newtheorem{re}[de]{Remark}

\newtheorem{te}[de]{Theorem}
\newtheorem{ex}[de]{Example}

\def\Box{\mbox{$\sqcap\!\!\!\!\sqcup$}}

\def\ot{\otimes}
\def\ra{\rightarrow}
\def\lra{\longrightarrow}
\def\al{\alpha}

\def\bea{\begin{eqnarray*}}
\def\eea{\end{eqnarray*}}

\def\tie{\mathrel>\joinrel\mathrel\triangleleft}
\def\Cut{\mbox{$\sqcap \!\!\!\! \sqcup$}}

\begin{document}
\title{Wide Morita contexts, relative injectivity and equivalence results}
\author{N. Chifan$^1$, S.
D\u{a}sc\u{a}lescu$^2$\thanks{ On leave from University of
Bucharest, Dept. Mathematics. } and
C. N\u{a}st\u{a}sescu$^3$\\[2mm]
$^1$ Facultatea de Matematic\u{a}, University of Gala\c{t}i, Romania\\
$^2$ Kuwait University, Faculty of Science, Dept. Mathematics,
PO BOX 5969,\\ Safat 13060, Kuwait, e-mail: sdascal@mcs.sci.kuniv.edu.kw\\
$^3$  University of Bucharest, Facultatea de Matematica, Str.
Academiei 14,\\ Bucharest 1, RO-010014, Romania, e-mail:
cnastase@al.math.unibuc.ro }

\date{}
\maketitle

\begin{abstract}
We extend Morita theory to abelian categories by using wide Morita
contexts. Several equivalence results are given for wide Morita
contexts between abelian categories, widely extending equivalence
theorems for categories of modules and comodules due to Kato,
M\"{u}ller and Berbec. In the case of Grothendieck categories we
derive equivalence results by using quotient categories. We apply
the general equivalence results to rings with identity, rings with
local units, graded rings, Doi-Hopf modules and coalgebras.\\[2mm]
\end{abstract}

\section{Introduction and preliminaries}

Morita contexts appeared in the work of Morita on equivalence of
categories of modules over rings with identity. A fundamental
result of Morita says that the categories of modules over two
rings with identity $R$ and $S$ are equivalent if and only if
there exists a strict Morita context connecting $R$ and $S$.
Morita contexts have been used to the study of group actions on
rings and Galois theory for commutative rings (see
\cite{meyeringraham}). A Morita theory for rings with local units
was developed in \cite{am}. Several Morita contexts were
constructed in connection to Galois theory for Hopf algebra
actions and coactions (see \cite{cfm}, \cite{cc}, \cite{bdr}),
where Hopf-Galois extensions are characterized by the surjectivity
of one of the Morita maps. As an application, the finite
dimensional version of the duality theorem of Blattner-Montgomery
(\cite{bm}) was deduced and explained in a nice way by using
Morita contexts and Hopf-Galois theory in \cite{vdb}.

A Morita context gives rise to an equivalence of categories if and
only if it is strict, i.e. if both Morita maps are surjective.  A
natural question that was posed was how far is an arbitrary Morita
context from an equivalence of categories. An answer is given by
the Kato-M\"{u}ller Theorem (see \cite{kato0}, \cite{muller}),
which briefly says (in the formulation of M\"{u}ller) that if the
rings $R$ and $S$ are connected by a Morita context, then certain
quotient categories of $R-mod$ and $S-mod$ are equivalent.

A concept dual to Morita contexts was constructed for coalgebras
by Takeuchi in \cite{ta}, where he defines what is now known as a
Morita-Takeuchi context connecting two coalgebras, and proves that
the categories of comodules over two coalgebras $C$ and $D$ are
equivalent if and only if $C$ and $D$ are connected by a strict
Morita-Takeuchi context. A result dual to the Kato-M\"{u}ller
Theorem for Morita-Takeuchi contexts and the associated categories
of comodules was proved by Berbec in \cite{berbec}.

Inspired by \cite{nw}, where an equivalence result for the
subcategories of $R-mod$ and $S-mod$ consisting of the
trace-torsionfree trace-accessible modules was proved in the case
where the rings $R$ and $S$ are connected by a Morita context,
Casta\~{n}o Iglesias and Gomez Torrecillas define in \cite{cg0},
\cite{cg} the concept of a wide Morita context for abelian
categories. A datum $(F,G,\eta ,\rho )$ is called a right wide
Morita context between the abelian categories $\cal A$ and $\cal
B$ if $F:{\cal B}\ra {\cal A}$ and $G:{\cal A}\ra {\cal B}$ are
right exact functors, and $\eta :F\circ G\ra 1_{\cal A}$ and $\rho
:G\circ F\ra 1_{\cal B}$ are natural transformations with the
property that $G\eta =\rho G$ and $F\rho =\eta F$. A left wide
Morita context is a datum which is a right wide Morita context
when regarded between the dual categories. Both Morita contexts
and Morita-Takeuchi contexts can be regarded as particular wide
Morita contexts, so this approach can be seen as a way to unify
(as much as possible) Morita theory for modules and comodules.

The aim of this paper is two-fold. On one hand, we study general
properties of wide Morita contexts and we explain how they allow
extending Morita theory to abelian categories. As a generalization
of Morita's result, we show that an equivalence of abelian
categories is essentially a strict wide Morita context. On the
other hand we give several equivalence results which widely extend
the results of Kato, M\"{u}ller and Berbec. We give a general
equivalence result for wide Morita contexts between abelian
categories, and we derive several equivalence theorems for the
case where the categories are Grothendieck. The M\"{u}ller type
equivalence result for Grothendieck categories seems to be the
most interesting since the quotient categories are also
Grothendieck categories, and this general framework can be applied
to a large number of examples. We apply the general results about
wide Morita contexts to rings with identity, rings with local
units, graded rings, Doi-Hopf modules and coalgebras. Some of the
results obtained in this way are known, but we explain them from a
general point of view, and some others are new.

The content of the paper is as follows. In Section 1 we recall the
definition of right wide Morita contexts and left wide Morita
contexts, and present several general properties of these. We also
define a composition for wide Morita contexts and show that it is
associative. We define the concepts of isomorphic right wide
Morita contexts, and we show that the invertible right wide Morita
contexts (with respect to the composition previously defined) are
exactly the strict right wide Morita contexts, i.e. those ones
defining an equivalence of categories. On the other hand, we show
that any equivalence of abelian categories arises from a strict
right wide Morita context. A fundamental fact that we prove is
that to a right wide Morita context for which $F$ and $G$ are left
adjoint functors, we can associate a left wide Morita context
whose functors are the right adjoints of $F$ and $G$.  In Section
2 we consider the concept of relative injective object with
respect to a subcategory and another object, and study related
properties. This allows us to define a concept of a closed object
with respect to a subcategory. For any right wide Morita context
$\Gamma$ between the abelian categories $\cal A$ and $\cal B$ such
that the functors $F$ and $G$ have right adjoints, we construct
subcategories ${\cal C}_{\Gamma}$ of $\cal A$, and ${\cal
D}_{\Gamma}$ of $\cal B$ such that the category of ${\cal
C}_{\Gamma}$-closed objects of $\cal A$ and the category of ${\cal
D}_{\Gamma}$-closed objects of $\cal B$ are equivalent. In Section
3 we consider the dual situation, by defining relative projective
objects and the dual results for left wide Morita contexts. Since
the dual of an abelian category is abelian, these results are dual
to the ones in Section 2, so they follow directly by dualization,
and we did not include direct proofs. Even if they follow directly
by dualization, we include them since they are interesting for
several applications. In Section 4 we consider that $\cal A$ and
$\cal B$ are Grothendieck categories. Then it makes sense to
consider the smallest localizing category $\overline{{\cal
C}_{\Gamma}}$ that contains ${\cal C}_{\Gamma}$, and discuss the
connection between $\overline{{\cal C}_{\Gamma}}$-closed objects
and ${\cal C}_{\Gamma}$-closed objects. We prove a general
equivalence result by using quotient categories and the results of
Section 2. In Section 5 we apply the general equivalence result to
some particular cases. If $\cal A$ and $\cal B$ are categories of
modules over rings with identity, we obtain as a particular case
Kato-M\"{u}ller  Theorem. We get new similar equivalence results
for Morita contexts associated to graded rings and to rings with
local units. The results of Section 3 are applied to left wide
Morita contexts obtained by taking the Hom functors associated to
a Morita context. We obtain results of Kato and Ohtake as
particular cases. We also apply the general equivalence result to
Doi-Hopf modules. As particular situations we obtain equivalence
results for Hopf-Galois extensions and for the case where a total
integral exists. We derive the Weak Structure Theorem for
Hopf-Galois extension as a special case. In Section 6 we consider
left wide Morita contexts between Grothendieck categories. In the
case the functors $F$ and $G$ commute with direct limits we obtain
a new equivalence result. In Section 7 we apply it to
Morita-Takeuchi contexts and derive as a particular case the
theorem of Berbec, and also to Hopf-Galois coextensions.

For notations and basic concepts we refer to \cite{macl} for
general category theory issues, to \cite{gabriel}, \cite{F} and
\cite{stenstrom} for things related to abelian categories,
Grothendieck categories and quotient categories, and to \cite{dnr}
and \cite{mo} for coalgebras, Hopf algebras and Hopf-Galois
theory. If $\cal A$ is a category, by a subcategory of $\cal A$ we
always mean a full subcategory. If $\cal A$ is an abelian
category, then the subcategory $\cal C$ of $\cal A$ is closed if
it is closed under subobjects, factor objects, and arbitrary
direct sums. If moreover $\cal C$ is closed under extensions, it
is called a localizing subcategory. By functor we always mean a
covariant functor. If $f:X\ra Y$ and $g:Y\ra Z$ are two morphisms
in a category, their composition is denoted by $g\circ f$. The
same notation is used for composition of functors. All the
categories we work with are abelian, and all functors are
additive. If $T,S:{\cal A}\ra {\cal B}$ are two functors, then a
natural transformation $\eta :T\ra S$ is called an epimorphism
(monomorphism) if $\eta (X)$ is an epimorphism (monomorphism) for
any $X\in {\cal A}$.

\section{Wide Morita contexts, equivalence of abelian categories
and adjunctions}

Let $\cal A$ and $\cal B$ be two abelian categories. Following
\cite{cg0}, \cite{cg} a datum $(F,G,\eta ,\rho )$ is called a
right wide Morita context between the categories $\cal A$ and
$\cal B$ if $F:{\cal B}\ra {\cal A}$ and $G:{\cal A}\ra {\cal B}$
are right exact functors, and $\eta :F\circ G\ra 1_{\cal A}$ and
$\rho :G\circ F\ra 1_{\cal B}$ are natural transformations with
the property that $G\eta =\rho G$ and $F\rho =\eta F$. Note that
in this case $(G,F,\rho ,\eta )$ is a right wide Morita context
between the categories $\cal B$ and $\cal A$, therefore any
general result that we prove for $F$, respectively
$\eta$, also holds for $G$, respectively $\rho$.\\
Dually, $(F,G,\eta ,\rho )$ is called a left wide Morita context
between $\cal A$ and $\cal B$ if $F:{\cal B}\ra {\cal A}$ and
$G:{\cal A}\ra {\cal B}$ are left exact functors, $\eta :1_{\cal
A}\ra F\circ G$ and $\rho :1_{\cal B}\ra G\circ F$ satisfy $G\eta
=\rho G$ and $F\rho =\eta F$.

\begin{pr} \label{etaepi}
Let $\Gamma =(F,G,\eta ,\rho )$ be a right wide Morita context. If
$\eta$ is an epimorphisms, then $\eta$ is a natural equivalence.
Thus if $\eta$ and $\rho$ are epimorphisms, the functors $F$ and
$G$ give an equivalence between the categories $\cal A$ and $\cal
B$.
\end{pr}
{\bf Proof:} Assume that $\eta$ is an epimorphism. Let $M\in {\cal
A}$. We have the exact sequence
$$
0\lra {\rm Ker}\; \eta (M)\stackrel{i}{\lra} (F\circ
G)(M)\stackrel{\eta (M)}{\lra} M\lra 0$$ where $i$ is the
inclusion morphism. Since $F\circ G$ is right exact, we obtain the
commutative diagram

\begin{picture}(100,100)(10,10)
\put(-7,80){$(F\circ G)({\rm Ker}\; \eta (M))$} \put(-2,47){$\eta
({\rm Ker}\; \eta (M))$} \put(32,10){${\rm Ker} \; \eta (M)$}

\put(8,10){$0$}\put(20,14){\vector(1,0){10}}

\put(65,75){\vector(0,-1){53}} \put(155,80){$(F\circ G)((F\circ
G)(M))$} \put(183,47){$\eta ((F\circ G)(M))$}
\put(155,10){$(F\circ G)(M)$} \put(177,75){\vector(0,-1){53}}
\put(87,14){\vector(1,0){65}} \put(95,84){\vector(1,0){57}}
\put(124,18){$i$} \put(95,88){$(F\circ G)(i)$}

\put(266,84){\vector(1,0){83}} \put(270,88){$(F\circ G)(\eta
(M))$}

\put(352,80){$(F\circ G)(M)$} \put(420,84){\vector(1,0){20}}
 \put(450,80){$0$}

\put(367,75){\vector(0,-1){53}} \put(372,47){$\eta (M)$}

\put(220,14){\vector(1,0){133}} \put(290,18){$\eta (M)$}

\put(362,10){$M$} \put(390,14){\vector(1,0){48}}
 \put(450,10){$0$}
\end{picture}

But \bea \eta ((F\circ G)(M))&=&F(\rho (G(M)))\\
&=&(F\circ G)(\eta (M))\eea and \bea \eta ((F\circ G)(M))\circ
(F\circ G)(i)&=&(F\circ G)(\eta (M))\circ (F\circ
G)(i)\\&=&(F\circ G)(\eta (M)\circ i)\\&=&(F\circ G)(0)\\&=&0\eea
so $i\circ \eta ({\rm Ker}\; \eta (M))=0$. Since $i$ is a
monomorphism, we have $\eta ({\rm Ker}\; \eta (M))=0$, and by the
hypothesis we get ${\rm Ker}\; \eta (M)=0$. We conclude that $\eta
(M)$ is an isomorphism. Similarly (or by using the remark that
general facts about $\eta$ also hold for $\rho$) one proves that
$\rho$ is an isomorphism whenever it is an epimorphism, so the
last part of the statement follows.\qed

The terminology of the following definition is inspired by the one
for classical Morita contexts. It will be clear in Section
\ref{sectionapplications} that in fact wide Morita contexts
generalize classical Morita contexts.

\begin{de}
A right wide Morita context $\Gamma =(F,G,\eta ,\rho )$ is called
strict if $\eta$ and $\rho$ are epimorphisms (so then by
Proposition \ref{etaepi}, $F$ and $G$ define an equivalence
between $\cal A$ and $\cal B$ with natural equivalences $\eta
:F\circ G\ra 1_{\cal A}$ and $\rho :G\circ F\ra 1_{\cal B}$).
\end{de}

\begin{re}
Let $\Gamma =(F,G,\eta ,\rho )$ be a right wide Morita context
between the categories $\cal A$ and $\cal B$. Then we can regard
$F$ and $G$ as functors between the dual categories ${\cal A}^0$
and ${\cal B}^0$. It is easy to see that in this way $\Gamma$
becomes a left wide Morita context between the categories ${\cal
A}^0$ and ${\cal B}^0$. Note that ${\cal A}^0$ and ${\cal B}^0$
are also abelian categories. \\
Similarly, any left wide Morita context can be regarded as a right
wide Morita context between the dual categories. In this way we
will be able to transfer results from right to left wide Morita
contexts.
\end{re}

A first example of how the above remark can be applied is the
following.

\begin{pr}
Let $(F,G,\eta ,\rho )$ be a left wide Morita context. If $\eta$
(respectively $\rho$) is a monomorphism, then it is a natural
equivalence.
\end{pr}

Now we define a composition operation for right wide Morita
contexts. Let $\Gamma =(F,G,\eta ,\rho )$ be a right wide Morita
context between the categories $\cal A$ and $\cal B$, and let
$\Delta =(U,V,\epsilon ,\delta )$ be a right wide Morita context
between the categories $\cal B$ and $\cal C$. Define the natural
transformations
$$\gamma :F\circ U\circ V\circ G\ra 1_{\cal A},\;\; \gamma
(X)=\eta (X)\circ F(\epsilon (G(X)))\;\; {\rm for}\; X\in {\cal
A}$$
$$\pi :V\circ G\circ F\circ U\ra 1_{\cal C},\;\; \pi
(Z)=\delta (Z)\circ V(\rho (U(Z)))\;\; {\rm for}\; Z\in {\cal C}$$
Shortly we denote $\gamma =\eta \circ F\epsilon G$ and $\pi
=\delta \circ V\rho U$. With these notations we have

\begin{pr} \label{prcompunere}
$(F\circ U,V\circ G,\gamma ,\pi )$ is a right wide Morita context
between the categories $\cal A$ and $\cal C$. We call this context
the composition of $\Gamma$ and $\Delta$, and we denote it by
$\Gamma \circ \Delta$.
\end{pr}
{\bf Proof:} We first note that since $\epsilon :U\circ V\ra
1_{\cal B}$ is a natural transformation, then for any morphism
$u:Y_1\ra Y_2$ in the category $\cal B$, we have
\begin{equation} \label{epsilonnatural}
\epsilon (Y_2)\circ (U\circ V)(u)=u\circ \epsilon (Y_1)
\end{equation}
Let $Z\in {\cal C}$. Then we have

\bea (F\circ U)(\pi (Z))&=&(F\circ U)(\delta (Z))\circ (F\circ
U\circ V)(\rho (U(Z)))\\&=&F(\epsilon (U(Z)))\circ F((U\circ
V)(\rho (U(Z))))\\
&=&F(\epsilon (U(Z))\circ (U\circ V)(\rho (U(Z))))\\
&=&F(\rho (U(Z))\circ \epsilon ((G\circ F\circ U)(Z)))\;\;\; ({\rm
by}\; \ref{epsilonnatural})\\
&=&\eta (F(U(Z)))\circ F(\epsilon ((G\circ F\circ U)(Z)))\\
&=&\gamma ((F\circ U)(Z))\eea showing that $(F\circ U)\pi =\gamma
(F\circ U)$. In a similar way we can prove that $(V\circ G)\gamma
=\pi (V\circ G)$. \qed

The composition of right wide Morita contexts is associative, as
the following result shows.

\begin{pr}
Let us consider three right wide Morita contexts: $\Gamma$ from
$\cal A$ to $\cal B$, $\Delta$ from $\cal B$ to $\cal C$, and
$\Sigma$ from $\cal C$ to $\cal D$. Then $(\Gamma \circ \Delta
)\circ \Sigma =\Gamma \circ (\Delta \circ \Sigma)$.
\end{pr}
{\bf Proof:} If $\Gamma =(F,G,\eta ,\rho )$, $\Delta
=(U,V,\epsilon ,\delta )$, and $\Sigma =(P,Q,\al ,\beta )$, it
follows directly from the definition that $$(\Gamma \circ \Delta
)\circ \Sigma =\Gamma \circ (\Delta \circ \Sigma)=(F\circ U\circ
P,Q\circ V\circ G, \eta \circ F\epsilon G\circ FU\al VG, \beta
\circ Q\delta P\circ QV\rho UP)$$ \qed

If $\cal A$ is an abelian category, $1_{\cal A}:{\cal A}\ra {\cal
A}$ is the identity functor, and $Id_{1_{\cal A}}:1_{\cal A}\ra
1_{\cal A}$ is the identity natural transformation, then clearly
we have a right (and also left) wide Morita context ${\bf 1}_{\cal
A}=(1_{\cal A}, 1_{\cal A}, Id_{1_{\cal A}}, Id_{1_{\cal A}})$
from $\cal A$ to $\cal A$. We call ${\bf 1}_{\cal A}$ the identity
wide Morita context. It is obvious that for any right wide Morita
context $\Gamma$ from $\cal A$ to $\cal B$ we have ${\bf 1}_{\cal
A}\circ \Gamma =\Gamma$ and $\Gamma \circ {\bf 1}_{\cal
B}=\Gamma$. Now we define a concept of isomorphism between wide
Morita contexts.

\begin{de} \label{defiso}
Let $\Gamma =(F,G,\eta ,\rho )$ and $\Delta =(P,Q,\al ,\beta )$ be
two right wide Morita contexts between $\cal A$ and $\cal B$. We
say that $\Gamma$ and $\Delta$ are isomorphic, and we write
$\Gamma \simeq \Delta$, if there exist natural equivalences
$u:F\ra P$ and $v:G\ra Q$ such that for any $X\in {\cal A}$ the
diagram

\begin{picture}(100,160)(-40,-20)
\put(70,120){$F(G(X))$} \put(70,54){$P(G(X))$}
\put(98,115){\vector(0,-1){49}} \put(101,90){$u(G(X))$}
 \put(257,120){$X$}
 \put(261,115){\vector(0,-1){108}}\put(264,54){$1_X$}
 \put(127,124){\vector(1,0){123}}
\put(166,128){$\eta (X)$} \put(98,49){\vector(0,-1){46}}
\put(101,25){$P(v(X))$} \put(70,-10){$P(Q(X))$}
\put(166,-2){$\al(X)$}\put(127,-6){\vector(1,0){123}}
\put(257,-10){$X$}
\end{picture}

is commutative, and for any $Y\in {\cal B}$ the diagram

\begin{picture}(100,160)(-40,-20)
\put(70,120){$G(F(Y))$} \put(70,54){$Q(F(Y))$}
\put(98,115){\vector(0,-1){49}} \put(101,90){$v(F(Y))$}
 \put(257,120){$Y$}
 \put(261,115){\vector(0,-1){108}}\put(264,54){$1_Y$}
 \put(127,124){\vector(1,0){123}}
\put(166,128){$\rho (Y)$} \put(98,49){\vector(0,-1){46}}
\put(101,25){$Q(u(Y))$} \put(70,-10){$Q(P(Y))$}
\put(166,-2){$\beta (Y)$}\put(127,-6){\vector(1,0){123}}
\put(257,-10){$Y$}
\end{picture}

is commutative. Shortly, we write these as $\al \circ Pv\circ
uG=\eta$ and $\beta \circ Qu\circ vF=\rho$.
\end{de}

\begin{re}  \label{remisoepi}
With the notations from Definition \ref{defiso} we have the
following.\\
(1) If $\Gamma \simeq \Delta$, then $\Delta \simeq \Gamma$, with
the latter isomorphism defined by the natural equivalences
$u^{-1}$ and $v^{-1}$.\\
(2) Assume that $\Gamma \simeq \Delta $. If $\eta$ is an
epimorphism, then so is $\al$. Also, if $\rho$ is an epimorphism,
then $\beta$ is an epimorphism.
\end{re}

\begin{de}
A right wide Morita context $\Gamma$ between $\cal A$ and $\cal B$
is called invertible if there exists a wide right Morita context
$\Delta$ between $\cal B$ and $\cal A$ such that $\Gamma \circ
\Delta \simeq {\bf 1}_{\cal A}$ and $\Delta \circ \Gamma\simeq
{\bf 1}_{\cal B}$.
\end{de}

The following result shows that the invertible right wide Morita
contexts are exactly the strict right wide Morita contexts.

\begin{pr}
Let $\Gamma =(F,G,\eta ,\rho )$ be a right wide Morita context
between the abelian categories $\cal A$ and $\cal B$. Then the
following assertions are equivalent.\\
(1) $\Gamma$ is invertible.\\
(2) $\Gamma$ is a strict right wide Morita context, i.e. $F$ and
$G$ define an equivalence between the categories $\cal A$ and
$\cal B$ with natural equivalences $\eta :F\circ G\ra 1_{\cal A}$
and $\rho :G\circ F\ra 1_{\cal B}$.
\end{pr}
{\bf Proof:} (1)$\Rightarrow$(2) Assume that there exists a right
wide Morita context $\Sigma =(U,V,\epsilon ,\delta )$ between
$\cal B$ and $\cal A$ such that $\Gamma \circ \Sigma \simeq {\bf
1}_{\cal A}$ and $\Sigma \circ \Gamma \simeq {\bf 1}_{\cal B}$.
Thus $(F\circ U,V\circ G,\eta \circ F\epsilon G ,\delta \circ
V\rho U)\simeq (1_{\cal A}, 1_{\cal A}, Id_{1_{\cal A}},
Id_{1_{\cal A}})$. By Remark \ref{remisoepi}(2) we have that $\eta
$ and $\delta$ are epimorphisms. Similarly we get that $\rho$ and
$\epsilon$ are epimorphisms, so $\Gamma$ and $\Sigma$ are strict.

(2)$\Rightarrow$(1) Let $\Delta =(G,F,\rho ,\eta)$, a right wide
Morita context. Then we have that
$$\Gamma \circ \Delta =(F\circ G,
F\circ G, \eta \circ F\rho G, \eta \circ F\rho G)\simeq {\bf
1}_{\cal A}$$ Indeed, we have the natural equivalences $\eta
:F\circ G\ra 1_{\cal A}$ on the role of $u$, and also $\eta$ on
the role of $v$ from Definition \ref{defiso}. Moreover
$$Id_{1_{\cal A}}\circ 1_{\cal A}\eta \circ \eta (F\circ G)=\eta
\circ \eta (F\circ G)=\eta \circ F\rho G$$ Similarly $\Delta \circ
\Gamma \simeq 1_{\cal B}$, and this ends the proof. \qed

We have seen that a strict right wide Morita context defines an
equivalence of categories. The next result shows that any
equivalence between abelian categories arises like this. These
generalize the classical result of Morita which tells that two
categories of modules are equivalent if and only if there exists a
strict Morita context connecting them.

\begin{pr}
Let $\cal A$ and $\cal B$ be abelian categories and $F:{\cal B}\ra
{\cal A}$, $G:{\cal A}\ra {\cal B}$ two functors defining an
equivalence of categories. Then there exists a wide Morita context
$(F,G,\eta ,\rho )$ between $\cal A$ and $\cal B$.
\end{pr}
{\bf Proof:} Since $F$ and $G$ define an equivalence of
categories, there exist two natural equivalences $u:F\circ G\ra
1_{\cal A}$ and $v:G\circ F\ra 1_{\cal B}$. Let $Y\in {cal B}$.
Then $u(F(Y)):(F\circ G\circ F)(Y)\ra F(Y)$ is a morphism in $\cal
A$. Since $F$ is fully faithful, there exists a unique morphism
$w(Y):(G\circ F)(Y)\ra Y$ in $\cal B$ such that $F(w(Y)=u(F(Y))$.
We show that $w:G\circ F\ra 1_{\cal B}$ is a natural
transformation, i.e. for any morphism $f:Y_1\ra Y_2$ in $\cal B$
we have $w(Y_2)\circ (G\circ F)(f)=f\circ w(Y_1)$. Since $F$ is
fully faithful, this is equivalent to $F(w(Y_2))\circ (F\circ
G\circ F)(f)=F(f)\circ F(w(Y_1))$. By the definition of $w$, this
is the same to $u(F(Y_2))\circ (F\circ G)(F(f))=F(f)\circ
u(F(Y_1))$, and this is true since $u:F\circ G\ra 1_{\cal A}$ is a
natural transformation.

Since $uF=Fw$, $u$ is a natural equivalence and $F$ is fully
faithful, we get that $w$ is also a natural equivalence.

It remains to show that $Gu=wG$, and then we have that $(F,G,u,w)$
is a right wide Morita context. Denote $H=F\circ G$. Since $u:H\ra
1_{\cal A}$ is a natural transformation, we have that $u(X)\circ
H(u(X))=u(X)\circ u(H(X))$ for any $X\in {\cal A}$. Since $u(X)$
is an isomorphism, this shows that $uH=Hu$, i.e. $(F\circ
G)(u(X))=u((F\circ G)(X))$. Since $uF=Fw$, this implies that
$F(G(u(X)))=F(w(G(X)))$. But $F$ is fully faithful, so then
$G(u(X))=w(G(X))$, showing that $Gu=wG$, and this ends the proof.
\qed

We end the section by discussing right wide Morita contexts for
which the two functors are left adjoints. Let $\Gamma =(F,G,\eta
,\rho )$ be a right wide Morita context between the abelian
categories $\cal A$ and $\cal B$. Assume that $(F,F',\al ,\beta )$
is an adjunction, i.e. $F'$ is a right adjoint of $F$, and $\al
:1_{\cal B}\ra F'\circ F$ and $\beta :F\circ F'\ra 1_{\cal A}$ are
the unit and the counit of the adjunction. This means that both
composites $F'\stackrel{\al F'}{\lra} F'\circ F\circ
F'\stackrel{F'\beta }{\lra}F'$ and $F\stackrel{F\al }{\lra} F\circ
F'\circ F\stackrel{\beta F}{\lra}F$ are identities. Also consider
an adjunction $(G,G',\gamma ,\delta )$, so $\gamma :1_{\cal A}\ra
G'\circ G$, $\delta :G\circ G'\ra 1_{\cal B}$, and both composites
$G'\stackrel{\gamma G'}{\lra} G'\circ G\circ G'\stackrel{G'\delta
}{\lra}G'$ and $G\stackrel{G\gamma }{\lra} G\circ G'\circ
G\stackrel{\delta G}{\lra}G$ are identities. It is known that the
functors $F'$, $G'$ are left exact. By \cite[page 101]{macl} we
can define the compositions of these adjunctions
$$(F,F',\al ,\beta )\circ (G,G',\gamma ,\delta
)=(F\circ G,G'\circ F',G'\al G\circ \gamma ,\beta \circ F\delta
F')$$ and
$$(G,G',\gamma ,\delta
)\circ (F,F',\al ,\beta )=(G\circ F,F'\circ G',F'\gamma F\circ \al
,\delta \circ G\beta G' )$$ We recall that the composition of the
adjunctions is associative (\cite[page 102]{macl}). Now we can
prove the following key result.

\begin{pr}  \label{oppositecontext}
Let $\Gamma =(F,G,\eta ,\rho )$ be a right wide Morita context
between the abelian categories $\cal A$ and $\cal B$ such that the
functors $F$ and $G$ have right adjoints $F'$ and $G'$. Then there
exist natural transformation $\eta ':1_{\cal A}\ra G'\circ F'$ and
$\rho ':1_{\cal B}\ra F'\circ G'$ induced by $\Gamma$ such that
$(F',G',\eta ',\rho ')$ is a left wide Morita context.
\end{pr}
{\bf Proof:} We keep the notations above for the adjunctions. By
\cite[Theorem 2, page 98]{macl}, there exists a unique natural
transformation $\eta ':1_{\cal A}\ra G'\circ F'$ (called the
conjugate of $\eta$) such that the diagram

\begin{picture}(100,100)(-60,10)
\put(20,80){$Hom(M,N)$} \put(-35,47){$Hom(\eta (M),1_N)$}
\put(-20,10){$Hom((F\circ G)(M),N)$}
\put(48,75){\vector(0,-1){53}} \put(170,80){$Hom(M,N)$}
\put(189,47){$Hom(1_M,\eta '(N))$} \put(160,10){$Hom(M,(G'\circ
F')(N))$} \put(186,75){\vector(0,-1){53}}
\put(94,14){\vector(1,0){62}} \put(90,84){\vector(1,0){75}}
\put(122,16){$\Phi$} \put(116,86){$=$}
\end{picture}

is commutative for any $M,N\in {\cal A}$, where $\Phi$ is the
isomorphism associated to the natural transformation $G'\al G\circ
\gamma$ or $\beta \circ F\delta F'$.

Similarly the conjugate $\rho ':1_{\cal B}\ra F'\circ G'$ of
$\rho$ is defined. By \cite[Theorem 2, page 102]{macl}, we have
that the conjugate of $G\eta =1_G\circ \eta$ is $\eta '\circ
1_{G'}=\eta 'G'$, and the conjugate of $\rho G$ is $G'\rho '$.
Since $G\eta =\rho G$ and the conjugate is unique, we have that
$\eta 'G'=G'\rho '$. Similarly we see that $\rho 'F'=F'\eta '$,
which shows that $(F',G',\eta ',\rho ')$ is a left wide Morita
context. \qed

\section{Relative injectivity and right wide Morita contexts}
\label{relinj}

Let $\cal A$ be an abelian category, and let $\cal C$ be a
subcategory of $\cal A$. If $X$ and $M$ are two objects of $\cal
A$, we say that $M$ is ${\cal C}$-$X$-injective if for any
monomorphism $i:X'\ra X$ in $\cal A$ such that $X/X'\in {\cal C}$,
and any morphism $f:X'\ra M$, there exists $g:X\ra M$ such that
$g\circ i=f$. \\
We denote by ${\cal I}^{-1}(M,{\cal C})=\{ X\in {\cal A}|\; M\;
{\rm is}\; {\cal C}-X-{\rm injective}\}$. If ${\cal
I}^{-1}(M,{\cal C})={\cal A}$, i.e. $M$ is ${\cal
C}$-$X$-injective for any $X\in {\cal A}$, we simply say that $M$
is $\cal C$-injective.

\begin{pr}  \label{proprXinj}
With the above notations, the following assertions hold.\\
(1) ${\cal I}^{-1}(M,{\cal C})$ is closed under factor objects.\\
(2) If $\cal C$ is closed under extensions, $X\in {\cal
I}^{-1}(M,{\cal C})$, and $X'$ is a subobject of $X$ such that
$X/X'\in {\cal C}$, then $X'\in {\cal I}^{-1}(M,{\cal C})$.\\
(3) If $\cal A$ is a Grothendieck category and $\cal C$ is closed
under subobjects, then ${\cal I}^{-1}(M,{\cal C})$ is closed under
direct sums.
\end{pr}
{\bf Proof:} (1) Assume that $M$ is ${\cal C}$-$X$-injective, and
let $Y$ be a subobject of $X$. We show that $M$ is ${\cal
C}$-$X/Y$-injective. Let $X'/Y\leq X/Y$ such that
$\frac{X/Y}{X'/Y}\simeq X/X'\in {\cal C}$. Denote by $p:X\ra X/Y$
and $p':X'\ra X'/Y$ the projection morphisms, and by $i:X'\ra X$
and $j:X'/Y\ra X/Y$ the inclusion morphisms such that $p\circ
i=j\circ p'$. Since $X\in {\cal I}^{-1}(M,{\cal C})$, we see that
there exists $h:X\ra M$ such that $h\circ i=f\circ p'$. But
$h(Y)=(h\circ i)(Y)=(f\circ p')(Y)=0$, so $h$ factorizes through
$X/Y$, i.e. there exists $g:X/Y\ra M$ with $g\circ p=h$. Then we
have $g\circ j\circ p'=g\circ p\circ i=h\circ i=f\circ p'$, and
since $p'$ is an epimorphism we have that $g\circ j=f$, which
shows that $M$ is ${\cal
C}$-$X/Y$-injective.\\
(2) Let $K$ be a subobject of $X'$ such that $X'/K\in {\cal C}$,
and let $f:K\ra M$ be a morphism. Denote by $j:K\ra X'$ and
$i:X'\ra X$ the inclusion morphisms. We have the exact sequence
$$0\longrightarrow X'/K\longrightarrow X/K\longrightarrow
X/X'\longrightarrow 0$$ and since $\cal C$ is closed under
extensions we obtain that $X/K\in {\cal C}$. Since $M$ is ${\cal
C}$-$X$-injective, there exists $g:X\ra M$ such that $g\circ
i\circ j=f$. Then $g\circ i:X'\ra M$ and $(g\circ i)\circ j=f$, so
$M$ is ${\cal
C}$-$X'$-injective.\\
(3) Let $(X_i)_{i\in I}$ be a family in ${\cal I}^{-1}(M,{\cal
C})$, and let $X=\oplus _{i\in I}X_i$. Let $K\leq X$ with $X/K\in
{\cal C}$, and let $h:K\ra M$ be a morphism. The family
$${\cal F}=\{ f:L\ra M|K\leq L\leq X\; {\rm and}\; f_{|K}=h\}$$
is inductive when ordered in the obvious way by inclusion, and
then by Zorn's Lemma it has a maximal element $\overline{h}:N\ra
M$. Let $N_i=N\cap X_i$. We have that $X_i/N_i=X_i/N\cap X_i\simeq
X_i+N/N\leq X/N$, and since $\cal C$ is closed under subobjects,
we also have $X_i/N_i\in {\cal C}$.

Assume that $N\neq X$. Then there is $i\in I$ such that $X_i$ is
not a subobject of $N$, and then $N_i$ is not a subobject of
$X_i$. Let $q:N_i\ra N$ be the inclusion morphism. Since $M$ is
${\cal C}$-$X_i$-injective, there exists $u:X_i\ra M$ extending
$\overline{h}\circ q$. But the restrictions of $u$ and
$\overline{h}$ to $N_i$ are equal, and then it is easy to see that
there exists $\overline{\overline{h}}:N+X_i\ra M$ extending both
$u$ and $\overline{h}$ (the argument is exactly as in
\cite[Proposition 1.13]{an}). Then $\overline{\overline{h}}\in
{\cal F}$ and $\overline{h}<\overline{\overline{h}}$, a
contradiction. We conclude that $N$ must be the whole of $X$, and
then $M$ is ${\cal C}$-$X$-injective.  \qed

Let $\cal A$ and $\cal B$ be abelian categories and let $\Gamma
=(F,G,\eta ,\rho )$ be a right wide Morita context between $\cal
A$ and $\cal B$. For any $M\in {\cal A}$ we write $I_M={\rm Im}\;
\eta (M)$, and for any $N\in {\cal B}$ we denote $J_N={\rm Im}\;
\rho (N)$. If $f:M\ra M'$ is a morphism in $\cal A$, the
commutative diagram

\begin{picture}(100,100)(-70,10)
\put(48,80){$FG(M)$} \put(72,47){$\eta (M)$} \put(60,10){$M$}
\put(68,75){\vector(0,-1){53}} \put(170,80){$FG(M')$}
\put(189,47){$\eta (M')$} \put(180,10){$M'$}
\put(186,75){\vector(0,-1){53}} \put(77,14){\vector(1,0){97}}
\put(90,84){\vector(1,0){75}} \put(122,18){$f$}
\put(116,88){$FG(f)$}
\end{picture}

shows that $f(I_M)\subseteq I_{M'}$. A similar result holds for
the objects $J_N, N\in {\cal B}$. Moreover, if $f:M\ra M'$ is an
epimorphism, we have $f(I_M)=I_{M'}$ since the functors $F$ and
$G$ are right exact. In particular we have that $I_{M/I_M}=0$.

We consider the following class \bea {\cal C}_{\Gamma}&=&\{ X\in
{\cal A}| f(I_M)=0 \; {\rm for}\; {\rm any}\; M\; {\rm and}\; f:M\ra X\}\\
&=&\{ X\in {\cal A}| f\circ \eta (M)=0 \; {\rm for}\; {\rm any}\;
M\; {\rm and}\; f:M\ra X\}\eea

\begin{pr} \label{propright}
With the above notations, the following assertions
hold.\\
(1) For any $M\in {\cal A}$ we have $M/I_M\in {\cal
C}_{\Gamma}$.\\
(2) $M\in {\cal C}_{\Gamma}$ if and only if $I_M=0$. Moreover, if
$M\in {\cal A}$ and $K\leq M$ is a subobject such that $M/K\in
{\cal C}_{\Gamma}$, then $I_M\leq K$.\\
(3) ${\cal C}_{\Gamma}$ is closed under subobjects and factor
objects. If $\cal A$ has direct products, then ${\cal C}_{\Gamma}$
is closed under direct products (i.e. ${\cal C}_{\Gamma}$ is a TTF-class).
Moreover, if $\cal A$ and $\cal B$ have direct sums and $F,G$ commute with
the direct sums, then ${\cal C}_{\Gamma}$ is closed under direct sums.\\
(4) If the category $\cal A$ has a family of generators
$(U_i)_{i\in I}$, then $\{ U_i/I_{U_i}|i\in I\}$ is a family of
generators
of the category ${\cal C}_{\Gamma}$.\\
(5) If $I_{I_M}=I_M$ for any $M$, then ${\cal C}_{\Gamma}$ is
closed under extensions.
\end{pr}
{\bf Proof:} (1) Let $N\in {\cal A}$ and let $f:N\ra M/I_M$ be a
morphism. Since $f(I_N)\leq I_{M/I_M}=0$, so $f(I_N)=0$, showing
that
$M/I_M\in {\cal C}_{\Gamma}$.\\
(2) If $I_M=0$, then we have $M\in {\cal C}_{\Gamma}$ by (1).
Conversely, let $M\in {\cal C}_{\Gamma}$. Then for the morphism
$1_M:M\ra M$ we have $0=1_M(I_M)=I_M$. The last part follows by
considering the natural projection $\pi :M\ra M/K$. Then $\pi
(I_M)=
I_{M/K}=0$, so $I_M\leq K$.\\
(3) Consider an exact sequence
$$0\lra X' \stackrel{i}{\lra} X\stackrel{\pi}{\lra} X"\lra 0$$
in ${\cal A}$, and assume that $X\in {\cal C}_{\Gamma}$. Let $M\in
{\cal A}$ and $f:M\ra X'$. Then $i\circ f:M\ra X$, so $(i\circ
f)(I_M)=0$. Since $i$ is a monomorphism, we have that $f(I_M)=0$,
so then $X'\in {\cal C}_{\Gamma}$.

Since $\pi$ is an epimorphism, we have that $\pi (I_X)=I_{X"}$.
But $X\in {\cal C}_{\Gamma}$ shows that $I_X=0$, and then we also
have $I_{X"}=0$, i.e. $X"\in {\cal C}_{\Gamma}$.

Let now $(X_i)_{i\in I}$ be a family of objects in ${\cal
C}_{\Gamma}$, $M\in {\cal A}$ and $f:M\ra \prod _{i\in I}X_i$ an
arbitrary morphism. If $\pi _i:\prod _{j\in I}X_j\ra X_i$ are the
natural projections, we have that $(\pi _i\circ f)(I_M)=0$ for any
$i$, so $f(I_M)=0$. We conclude that $\prod _{i\in I}X_i\in {\cal
C}_{\Gamma}$.

Assume now that $\cal A$ and $\cal B$ have direct sums, and $F$
and $G$ commute with direct sums. Let $(M_i)_{i\in I}$ be a family
of objects in ${\cal C}_{\Gamma}$, $\oplus _{i\in I}M_i$ their
direct sum, and $q_j:M_j\ra \oplus _{i\in I}M_i$, the natural
embedding for any $j\in I$. Since $F$ and $G$ commute with the
direct sums, so does $F\circ G$, therefore $(F\circ G)(\oplus
_{i\in I}M_i)\simeq \oplus _{i\in I}(F\circ G)(M_i)$. Since $\eta$
is a natural transformation we have that
$$\eta (\oplus _{i\in I}M_i)\circ (F\circ G)(q_j)=q_j\circ \eta
(M_j)=0$$ for any $j\in I$, and we conclude that $\eta (\oplus
_{i\in I}M_i)=0$, which shows that $\oplus _{i\in I}M_i\in {\cal
C}_{\Gamma}$.\\

(4) Let $X\in {\cal C}_{\Gamma}$ and $X'<X$ a strict subobject.
Since $(U_i)_{i\in I}$ is a family of generators of $\cal A$,
there exist $i\in I$ and a morphism $f:U_i\ra X$ such that ${\rm
Im}\; f$ is not a subobject of $X'$. Since $f(I_{U_i})\leq I_X=0$
(by (3)), then $f$ factorizes through a morphism $g:U_i/I_{U_i}\ra
X$. Clearly ${\rm Im}\; g={\rm Im}\; f$ is not a subobject of
$X'$, and this ends
the proof.\\
(5) Let
$$0\lra X' \stackrel{i}{\lra} X\stackrel{\pi}{\lra} X"\lra 0$$
be an exact sequence with $X',X"\in{\cal C}_{\Gamma}$. If $f:M\ra
X$ is a morphism, then $(\pi \circ f)(I_M)=0$, so there exists
$g:I_M\ra X'$ such that $i\circ g=f$. By hypothesis we have
$g(I_{I_M})=0$, so then $g(I_M)=0$. We conclude that $f(I_M)=0$
and $X\in {\cal C}_{\Gamma}$. \qed

\begin{re}
We note that part (1) of Proposition \ref{propright} shows that
${\cal C}_{\Gamma}=\{ X\in {\cal A}| \eta (X)=0\}$.
\end{re}

\begin{pr}  \label{kercgamma}
Let $M\in {\cal A}$. Then ${\rm Ker} \; \eta (M)\in {\cal
C}_{\Gamma}$.
\end{pr}
{\bf Proof:} Denote by $K={\rm Ker} \; \eta (M)$ and by $i:K\ra
(F\circ G)(M)$ the inclusion morphism. Let $U\in {\cal A}$ be an
object and $f:U\ra K$ a morphism. We have the commutative diagram

\begin{picture}(100,160)(-40,-20)
\put(70,120){$(F\circ G)(U)$} \put(70,54){$(F\circ G)(K)$}
\put(98,115){\vector(0,-1){49}} \put(101,90){$(F\circ G)(f)$}
 \put(257,120){$U$}
\put(257,54){$K$}
\put(261,115){\vector(0,-1){49}}\put(264,90){$f$}
\put(135,58){\vector(1,0){115}} \put(135,124){\vector(1,0){115}}
\put(166,62){$\eta (K)$} \put(166,128){$\eta (U)$}
\put(98,49){\vector(0,-1){49}} \put(101,25){$(F\circ G)(i)$}
\put(261,49){\vector(0,-1){49}}\put(264,25){$i$}
\put(35,-10){$(F\circ G)((F\circ G)(M))$} \put(150,-2){$\eta
((F\circ G)(M))$}\put(148,-6){\vector(1,0){87}}
\put(238,-10){$(F\circ G)(M)$}
\end{picture}

Since $\eta F=F\rho$ and $\rho G=G\eta$, we have $\eta ((F\circ
G)(M))=(F\circ G)(\eta (M))$, and then \bea i\circ f\circ \eta
(U)&=&i\circ \eta (K)\circ (F\circ G)(f)\\
&=&(F\circ G)(\eta (M))\circ (F\circ G)(i)\circ (F\circ G)(f)\\
&=&(F\circ G)(\eta (M)\circ i\circ f)\\
&=&0 \eea Since $i$ is a monomorphism we have $f\circ \eta (U)=0$,
showing that $K\in {\cal C}_{\Gamma}$. \qed

\begin{co}  \label{cogen}
If $(U_i)_{i\in I}$ is a family of generators of the category
$\cal A$ and $\eta (U_i)$ is an epimorphism for any $i\in I$, then
$\eta$ is a natural equivalence.
\end{co}
{\bf Proof:} By Proposition \ref{propright}(4), we see that ${\cal
C}_{\Gamma}=0$. Then by Proposition \ref{kercgamma} we have ${\rm
Ker}\; \eta (M)=0$ for any $M\in {\cal A}$, so $\eta$ is a natural
equivalence.\qed

\begin{de}
Let $M$ be an object of the abelian category $\cal A$, and $\cal
C$ be an arbitrary subcategory of $\cal A$. We say that $M$ is
$\cal C$-torsion free if for any $X\in \cal C$ and any
monomorphism $i:X\ra M$, we must have $X=0$. If $M$ is $\cal
C$-torsion free and $\cal C$-injective, then $M$ is called $\cal
C$-closed.
\end{de}

Now we can characterize the ${\cal C}_{\Gamma}$-closed objects by
a categorial property.

\begin{pr} \label{closed}
An object $M\in {\cal A}$ is ${\cal C}_{\Gamma}$-closed if and
only if for any $U\in {\cal A}$ the natural map
$$\phi =Hom(\eta (U),1_M):Hom(U,M)\ra Hom((F\circ G)(U),M)
,\;\; \phi (\beta )=\beta \circ \eta (U)$$ is bijective.
\end{pr}
{\bf Proof:} For $U\in {\cal A}$ we denote by $\eta_1 (U):(F\circ
G)(U)\ra I_U$ the corestriction of $\eta (U):(F\circ G)(U)\ra U$
to $I_U$, and by $j:I_U\ra U$ the inclusion morphism. Note that
$\eta (U)=j\circ \eta_1 (U)$.

Assume that $M$ is ${\cal C}_{\Gamma}$-closed. Let $\al :(F\circ
G)(U)\ra M$ be a morphism. We consider the following diagram.

\begin{picture}(100,110)(-30,-10)
\put(20,47){${\rm Ker}\; \eta (U)$}
\put(70,51){\vector(1,0){28}}\put(82,55){$i$}
 \put(100,47){$(F\circ G)(U)$}
\put(160,51){\vector(1,0){50}}\put(171,55){$\eta_1(U)$}
\put(215,47){$I_U$}\put(230,51){\vector(1,0){30}}
\put(241,55){$j$}\put(265,47){$U$} \put(128,43){\vector(0,-1){30}}
\put(122,0){$M$} \put(210,43){\vector(-2,-1){70}}
\put(260,43){\vector(-3,-1){120}} \put(120,23){$\al$}
\put(160,25){$\overline{\al}$} \put(200,10){$\beta$}
\end{picture}

As in the proof of Proposition \ref{kercgamma}, $i$ denotes the
inclusion morphism. Since ${\rm Ker}\; \eta (U)\in {\cal
C}_{\Gamma}$ (by Proposition \ref{kercgamma}) and $M$ is ${\cal
C}_{\Gamma}$-torsion free, we have that $\al ({\rm Ker}\; \eta
(U))=0$. Thus there exists a unique morphism
$\overline{\al}:I_U\ra M$ such that $\overline{\al}\circ \eta_1
(U)=\al$. Since $U/I_U\in {\cal C}_{\Gamma}$ and $M$ is ${\cal
C}_{\Gamma}$-injective, we see that there exists a morphism $\beta
:U\ra M$ such that $\overline{\al} =\beta \circ j$. Then $\al
=\beta \circ j\circ \eta_1(U)=\beta \circ \eta (U)=\phi (\beta )$,
showing that $\phi$ is surjective.

Now if $\phi (\beta)=\beta \circ \eta (U)=0$ for some $\beta$, we
have that $\beta (I_U)=0$, so then there exists
$\overline{\beta}:U/I_U\ra M$ such that $\overline{\beta}\circ
p=\beta$, where $p:U\ra U/I_U$ is the natural projection. But
$U/I_U\in {\cal C}_{\Gamma}$ and $M$ is ${\cal
C}_{\Gamma}$-torsion free, so we must have $\overline{\beta}=0$.
Hence $\beta =0$, so $\phi$ is injective.

Therefore $\phi$ a bijection, and it is obviously natural.

Conversely, assume that $\phi$ is bijective for any $U$. Define
$p:Hom(I_U,M)\ra Hom((F\circ G)(U),M)$ by $p(f)=f\circ \eta_1
(U)$, and $q:Hom(U,M)\ra Hom(I_U,M)$ by $q(f)=f\circ j$. Since
$\eta_1 (U)$ is an epimorphism, $p$ is injective. We have the
commutative diagram

\begin{picture}(100,120)(-20,20)
\put(45,90){$0$} \put(55,94){\vector(1,0){25}}
\put(184,90){$Hom((F\circ G)(U),M)$} \put(85,90){$Hom(I_U,M)$}
\put(163,98){$p$} \put(153,94){\vector(1,0){27}}
\put(105,50){\vector(0,1){35}} \put(96,66){$q$}
\put(88,35){$Hom(U,M)$} \put(150,42){\vector(2,1){68}}
\put(178,63){$\phi$}
\end{picture}

This shows that $p$ is surjective, so then it is bijective. Hence
$q$ is also bijective.

If $U$ is a subobject of $M$ which is in ${\cal C}_{\Gamma}$, we
have that $I_U=0$, and this shows that $Hom(U,M)=0$. Thus $M$ is
${\cal C}_{\Gamma}$-torsion free.

Now if we take an arbitrary $U$, and $K$ a subobject of $U$ such
that $U/K\in {\cal C}_{\Gamma}$, we have by Proposition
\ref{propright} that $I_U\leq K$. We thus have $I_K\leq I_U\leq
K\leq U$. If $f:K\ra M$ is a morphism, it restricts to a morphism
$f_1:I_U\ra M$. Since $q:Hom(U,M)\ra Hom(I_U,M)$ is an
isomorphism, there exists $h:U\ra M$ extending $f_1$. Denote by
$f'$ the restriction of $h$ to $K$. Let $f_2$ be the restriction
of $f$ to $I_K$ (which is also the restriction of $f_1$ to $I_K$).
Then clearly the restriction of $f'$ to $I_K$ is $f_2$, and the
natural isomorphism $Hom(I_K,M)\simeq Hom(K,M)$ shows that $f'=f$.
Thus $h$ extends $f$, and this proves that $M$ is ${\cal
C}_{\Gamma}$-injective.  \qed

We can prove now the main result of this section.

\begin{te}  \label{teinj}
Let $\Gamma =(F,G,\eta ,\rho )$ be a right wide Morita context
between the abelian categories $\cal A$ and $\cal B$, such that
the functors $F$ and $G$ have right adjoints $F'$ and $G'$. Denote
by ${\cal C}_{\Gamma}$ (respectively ${\cal D}_{\Gamma}$) the
subcategory of $\cal A$ (respectively $\cal B$) defined by the
natural morphism $\eta$ (respectively $\rho$). Then the categories
of ${\cal C}_{\Gamma}$-closed objects of $\cal A$ and ${\cal
D}_{\Gamma}$-closed objects of $\cal B$ are equivalent via the
functors $F',G'$.
\end{te}
{\bf Proof:} By Proposition \ref{oppositecontext}, we can
associate a left wide Morita context $(F',G',\eta ',\rho ')$ to
$\Gamma$. Let $M\in {\cal A}$ be ${\cal C}_{\Gamma}$-closed. By
Proposition \ref{closed} there is a natural bijection
$$\phi =Hom(\eta (U),1_M):Hom(U,M)\ra Hom((F\circ G)(U),M)
,\;\; \phi (\beta )=\beta \circ \eta (U)$$ where $U\in {\cal A}$.
Since $G'\circ F'$ is a right adjoint of $F\circ G$, we have a
natural bijection
$$\psi :Hom((F\circ G)(U),M)\ra Hom(U,(G'\circ F')(M))$$
The construction of $\eta '$ as the conjugate of $\eta$ (see the
proof of Proposition \ref{oppositecontext}) shows that $\psi \phi
=Hom(1_U, \eta '(M))$ Since $\psi$ and $\phi$ are bijections, then
so is $Hom(1_U,\eta '(M))$. This implies that $\eta '(M)$ is an
isomorphism. Similarly, if $N\in {\cal B}$ is ${\cal
D}_{\Gamma}$-closed, then the natural isomorphism $\rho '(N):N\ra
(F'\circ G')(N)$ is an isomorphism.

On the other hand, for any object $V\in {\cal B}$ we have a
commutative diagram

\begin{picture}(100,100)(-60,10)
\put(-20,80){$Hom(V,F'(M))$} \put(3,47){$\simeq$}
\put(-20,10){$Hom(F(V),M)$} \put(18,75){\vector(0,-1){53}}
\put(190,80){$Hom((G\circ F)(V),F'(M))$} \put(236,47){$\simeq$}
\put(190,10){$Hom((F\circ G\circ F)(V),M)$}
\put(231,75){\vector(0,-1){53}} \put(64,14){\vector(1,0){120}}
\put(64,84){\vector(1,0){120}} \put(81,19){$Hom(F(\rho (V)),1_M)$}
\put(81,89){$Hom(\rho (V), 1_{F'(M)})$}
\end{picture}

where the vertical arrows are the natural bijections coming from
the adjunction given by $F$ and $F'$. Since $F\rho =\eta F$ and
$M\in {\cal A}$ is ${\cal C}_{\Gamma}$-closed, we have that
$Hom(F(\rho (V)),1_M)=Hom(\eta (F(V)), 1_M)$ is a bijection, and
then also the top horizontal arrow $Hom(\rho (V),1_{F'(M)})$ is a
bijection. This shows that $F'(M)$ is ${\cal D}_{\Gamma}$-closed,
by Proposition \ref{closed}. Similarly $G'(N)$ is ${\cal
C}_{\Gamma}$-closed for any $N\in {\cal B}$ which is ${\cal
D}_{\Gamma}$-closed. These show that $F'$ and $G'$ induce an
equivalence between the subcategory of ${\cal C}_{\Gamma}$-closed
objects of $\cal A$ and the subcategory of ${\cal
D}_{\Gamma}$-closed objects of $\cal B$. \qed

\begin{co}  \label{corinj}
Let $\Gamma =(F,G,\eta ,\rho )$ be a right wide Morita context
between the abelian categories $\cal A$ and $\cal B$, such that
the functors $F$ and $G$ have right adjoints $F'$ and $G'$, and
$\eta$ is an epimorphism. Then $\cal A$ and the category of ${\cal
D}_{\Gamma}$-closed objects of $\cal B$ are equivalent via the
functors $F'$ and the restriction of $F$.
\end{co}
{\bf Proof:} Since $\eta$ is an isomorphism, we have ${\cal
C}_{\Gamma}=0$, so the category of all ${\cal C}_{\Gamma}$-closed
objects of $\cal A$ is the whole of $\cal A$. Denote by $\cal D$
the subcategory of $\cal B$ consisting of all ${\cal
D}_{\Gamma}$-closed objects. We know from Theorem \ref{teinj} that
$F'$ is an equivalence between $\cal A$ and $\cal D$. Since $F$ is
a left adjoint of $F'$, then the restriction $F:{\cal D}\ra {\cal
A}$ is still a left adjoint of $F'$, when regarded as a functor
from $\cal A$ to $\cal D$. But then $F'$ is an equivalence, so
$F:{\cal D}\ra {\cal A}$ is also an equivalence. \qed

\section{The dual case: relative projectivity and left wide Morita
contexts}  \label{relproj}

In this section we consider the dual concepts of the ones
presented in the Section \ref{relinj}. As we will see, this is
useful for understanding several examples. Since the dual of an
abelian category is also an abelian category, we can dualize the
definitions and results directly.

So let $\cal A$ be an abelian category and $\cal C$ a subcategory
of $\cal A$. If $M$ and $X$ are two objects of $\cal A$, we say
that $M$ is ${\cal C}$-$X$-projective if $M$ is ${\cal
C}^0$-$X$-injective as an object of the dual category ${\cal
A}^0$. It is clear that $M$ is ${\cal C}$-$X$-projective if for
any epimorphism $p:X\ra X'$ in $\cal A$, such that ${\rm Ker}\;
p\in {\cal C}$, and any morphism $f:P\ra X'$, there exists a
morphism $g:M\ra X$ such that $p\circ g=f$. We denote ${\cal
P}^{-1}(M,{\cal C})=\{ X\in {\cal A}|\; M\; {\rm is}\; {\cal
C}{\rm -}X{\rm -projective}\}$. If ${\cal P}^{-1}(M,{\cal
C})={\cal A}$, we say that $M$ is $\cal C$-projective. We have the
following properties.

\begin{pr} \label{proprXproj}
With the above notations, the following assertions hold.\\
(1) ${\cal P}^{-1}(M,{\cal C})$ is closed under subobjects. \\
(2) If $\cal C$ is closed under extensions, $X\in {\cal
P}^{-1}(M,{\cal C})$, and $p:X\ra X'$ is an epimorphism such that
${\rm Ker}\; p\in {\cal C}$, then $X'\in {\cal P}^{-1}(M,{\cal C})$.\\
(3) If $\cal A$ is a Grothendieck category and $\cal C$ is closed
under factor objects, then ${\cal P}^{-1}(M,{\cal C})$ is closed
under direct products. In particular ${\cal P}^{-1}(M,{\cal C})$
is
closed at finite direct sums.\\
(4) If $\cal A$ is a Grothendieck category, $\cal C$ is a closed
subcategory, and $M$ is a finitely generated object, then ${\cal
P}^{-1}(M,{\cal C})$ is closed under arbitrary direct sums.
\end{pr}
{\bf Proof:} Parts (1)-(3) follow directly by dualizing
Proposition \ref{proprXinj}. Part (4) can be proved as
\cite[Proposition 16.12]{af}. \qed

Let $M$ be an object of $\cal A$, and $\cal C$ be an arbitrary
subcategory of $\cal A$. We say that $M$ is $\cal C$-cotorsion
free if for any $X\in \cal C$ and any epimorphism $f:M\ra X$, we
have $f=0$. This is of course equivalent to the fact that $M$ is
${\cal C}^0$-torsion free when regarded in the dual category
${\cal A}^0$.

Let now $\cal A$ and $\cal B$ be abelian categories and let
$\Gamma =(F,G,\eta ,\rho )$ be a left wide Morita context between
$\cal A$ and $\cal B$. For any $M\in {\cal A}$ we denote $K_M={\rm
Ker} \; \eta (M)$, and for any $N\in {\cal B}$ we denote $L_N={\rm
Ker}\; \rho (N)$. By looking at the definition of $I_M$ and $J_N$
in Section \ref{relinj}, it is natural to consider these objects
$K_M$ and $L_N$, since when regarded in the dual category, the
image of a morphism becomes a coimage. If $f:M\ra M'$ is a
morphism in $\cal A$, we have that $f(K_M)\subseteq K_{M'}$.
Dually to the definition we made in Section \ref{relinj}, define
now \bea {\cal C}_{\Gamma}&=&\{ X\in {\cal A}|\eta (M)\circ f=0\;
{\rm for}\;
{\rm any}\; M\; {\rm and}\; f:X\ra M\}\\
&=&\{ X\in {\cal A}| {\rm Im}\;f\leq K_M \; {\rm for}\; {\rm
any}\; f:M\ra X\}\eea

If we regard $\Gamma$ as a right wide Morita context between the
dual categories, then ${\cal C}_{\Gamma}$ is exactly the class
associated to $\eta$. Similarly one can consider the subcategory
${\cal D}_{\Gamma}$ associated to $\rho$. The following results
are dual to Propositions \ref{propright} and \ref{kercgamma},
Corollary \ref{cogen}, and Proposition \ref{closed}.

\begin{pr} \label{propleft}
With the above notations, the following assertions
hold.\\
(1) For any $M\in {\cal A}$ we have $K_M\in {\cal
C}_{\Gamma}$.\\
(2) $M\in {\cal C}_{\Gamma}$ if and only if $K_M=M$.\\
(3) ${\cal C}_{\Gamma}$ is closed under subobjects and factor
objects. If $\cal A$ has direct sums, then ${\cal C}_{\Gamma}$ is
closed under direct sums (i.e. ${\cal C}_{\Gamma}$ is a closed
subcategory). Moreover, if $\cal A$ and $\cal B$ have direct
products and the functors $F$ and $G$ commute with direct
products, then ${\cal C}_{\Gamma}$ is closed under direct products
(i.e. it is a TTF-class). \\
(4) If the category $\cal A$ has a family of cogenerators
$(Q_i)_{i\in I}$, then $\{ K_{Q_i}|i\in I\}$ is a family of
cogenerators of the category ${\cal C}_{\Gamma}$.
\end{pr}

\begin{pr}  \label{cokercgamma}
Let $M\in {\cal A}$. Then ${\rm Coker} \; \eta (M)\in {\cal
C}_{\Gamma}$.
\end{pr}

\begin{co} \label{cocogen}
If $(Q_i)_{i\in I}$ is a family of cogenerators of the category
$\cal A$ and $\eta (Q_i)$ is a monomorphism for any $i\in I$, then
$\eta$ is a natural equivalence.
\end{co}

\begin{pr} \label{prclosed}
An object $M\in {\cal A}$ is ${\cal C}_{\Gamma}$-cotorsion free
and ${\cal C}_{\Gamma}$-projective if and only if for any $U\in
{\cal A}$, the map
$$\phi :Hom(M,U)\ra Hom(M,(F\circ G)(U))\;\; \phi (\beta )=\eta (U) \circ
\beta $$ is bijective.
\end{pr}

Denote by ${\cal A}_{\Gamma ,proj}$ the subcategory of $\cal A$
consisting of all objects that are ${\cal C}_{\Gamma}$-cotorsion
free and ${\cal C}_{\Gamma}$-projective.  Similarly we denote by
${\cal B}_{\Gamma ,proj}$ the subcategory of $\cal B$ consisting
of all objects that are ${\cal D}_{\Gamma}$-cotorsion free and
${\cal D}_{\Gamma}$-projective. The dual of Theorem \ref{teinj}
and Corollary \ref{corinj} are the following.

\begin{te} \label{teproj}
Let $\Gamma =(F,G,\eta ,\rho )$ be a left wide Morita context
between the abelian categories $\cal A$ and $\cal B$, such that
the functors $F$ and $G$ have left adjoints $F'$ and $G'$.  Then
the categories ${\cal A}_{\Gamma ,proj}$ and ${\cal B}_{\Gamma
,proj}$ are equivalent via the functors $F',G'$.
\end{te}

\begin{co}
Let $\Gamma =(F,G,\eta ,\rho )$ be a left wide Morita context
between the abelian categories $\cal A$ and $\cal B$, such that
the functors $F$ and $G$ have left adjoints $F'$ and $G'$. If
$\eta$ is a monomorphism, then the categories $\cal A$ and ${\cal
B}_{\Gamma ,proj}$ are equivalent via the functors $F'$ and the
restriction of $F$.
\end{co}

\section{Wide Morita contexts over Grothendieck categories and
equivalence results}

Throughout this section we assume that $\cal A$ is a Grothendieck
category and $\cal C$ is a closed subcategory of $\cal A$. Since
$\cal C$ is closed under factor objects, an object $M\in {\cal A}$
is $\cal C$-torsion free if for any object $X\in {\cal C}$ and any
morphism $f:X\ra M$, we have $f=0$. We denote by $t(M)$ the sum of
all subobjects of $M$ belonging to $\cal C$. Clearly $t(M)$ exists
since $\cal C$ is closed under arbitrary direct sums and factor
objects. In this way a left exact functor $t:{\cal A}\ra {\cal A}$
is defined; it is called the preradical associated to $\cal C$.
Clearly $M$ is $\cal C$-torsion free if and only if $t(M)=0$. If
$\cal C$ is a localizing subcategory, we have that $t(M/t(M))=0$.

For a closed subcategory $\cal C$, we denote by $\overline{\cal
C}$ the smallest localizing subcategory containing $\cal C$. This
is given by
$$\overline{\cal C}=\{ X\in {\cal A}|\; {\rm For}\; {\rm any}\; X'<X,
X/X'\; {\rm contains}\; {\rm a}\; {\rm non-zero}\; {\rm object}\;
{\rm of}\; {\cal C}\}$$ If $t$ is the preradical associated to
$\cal C$, we denote by $\overline{t}$ the preradical (which is in
fact a radical) associated to $\overline{\cal C}$. We have the
following characterization.

\begin{te}  \label{caractclosed}
Let $\cal C$ be a closed subcategory of $\cal A$, and let $M$ be
an object of $\cal A$. Then the following assertions are
equivalent.\\
(1) $M$ is $\overline{\cal C}$-closed.\\
(2) $M$ is $\cal C$-closed.\\
(3) If $(U_i)_{i\in I}$ is a family of generators of the category
$\cal A$, then $M$ is ${\cal C}$-$U_i$-injective for any $i\in I$,
and $t(M)=0$.
\end{te}
{\bf Proof:} (1)$\Rightarrow$(2) is clear since ${\cal C}\subseteq
\overline{\cal C}$.\\
(2)$\Rightarrow$(1) If $\overline{t}(M)\neq 0$, then by the
construction of $\overline{\cal C}$, the object $\overline{t}(M)$
contains a non-zero object belonging to $\cal C$, and this would
imply $t(M)\neq 0$, a contradiction. Therefore
$\overline{t}(M)=0$.

We prove now that $M$ is $\overline{\cal C}$-injective. Let $X\in
{\cal A}$ and $X'\leq X$ such that $X/X'\in \overline{\cal C}$,
and take $f:X'\ra M$ be a morphism. A standard application of
Zorn's Lemma shows that there exists a maximal subobject $Y$ of
$X$, with $X'\leq Y$ and there exists a morphism $g:Y\ra M$
extending $f$. Hence $X/Y$ is a factor object of $X/X'$, so
$X/Y\in \overline{\cal C}$, and then there exists a subobject
$Y<Z\leq X$ such that $Z/Y\in {\cal C}$. Since $M$ is $\cal
C$-injective, there exists $h:Z\ra M$ extending $g$. This is a
contradiction to the maximality of $Y$. We conclude that $Y=X$ and
$M$ is $\overline{\cal C}$-injective. \\
(2)$\Rightarrow$ (3) is clear.\\
(3)$\Rightarrow$(2) Let $X\in {\cal A}$. Since $X$ is a factor
object of a direct sum of $U_i$'s, and $M$ is ${\cal
C}$-$U_i$-injective for any $i$, we see by Proposition
\ref{proprXinj} that $M$ is ${\cal C}$-$X$-injective.  \qed

\begin{ex}  \label{exemple}
(i) Let $R$ be a ring with identity and $I$ be a two-sided ideal
of $R$. We define the class ${\cal P}_I=\{ M\in R-mod|IM=0\}$. It
is easy to see that ${\cal P}_I$ is a closed subcategory of
$R-mod$. The smallest localizing subcategory of $R-mod$ containing
${\cal P}_I$ is \bea {\cal C}_I&=&\{ M\in R-mod|\; {\rm For}\;
{\rm any}\; M'<M, M/M' \; {\rm contains}\; {\rm some}\; S\in {\cal
P}_I, S\neq 0\}\\
&=&\{ M\in R-mod|\; {\rm For}\; {\rm any}\; M'<M,\; {\rm there}\;
{\rm is}\; m\in M-M' \; {\rm such}\; {\rm that}\; Im\subseteq M'\}
\eea
Let us note that if the ideal $I$ is idempotent, i.e.
$I^2=I$, then ${\cal C}_I={\cal P}_I$. An object $M\in R-mod$ is
${\cal C}_I$-torsion free if and only if it is ${\cal
P}_I$-torsion free, and this means that
$$Ann_M(I)=\{ x\in
M|Ix=0\}=0$$
Now by Theorem \ref{caractclosed}, part (3), we have
that $M$ is ${\cal C}_I$-closed if and only if $Ann_M(I)=0$ and
any morphism $f:I\ra M$ of $R$-modules can be (uniquely) extended
to a morphism $g:R\ra M$. We conclude that $M$ is ${\cal
C}_I$-closed if and only if the natural morphism
$$\al :M\ra Hom_R(I,M),\;\; \al (m)(a)=am,\; a\in I, m\in M$$
is an isomorphism.

(ii) Let $R=\oplus _{\sigma \in G}R_{\sigma}$ be a $G$-graded
ring, and let $R-gr$ be the category of left graded $R$-modules.
If $M=\oplus _{\lambda \in G}M_{\lambda}$ is an object of this
category, we can consider for any $\sigma \in G$ the graded module
$M(\sigma )$ such that $M(\sigma )=M$ as an $R$-module, and the
homogeneous components are given by $M(\sigma
)_{\lambda}=M_{\lambda \sigma}$ for any $\lambda \in G$. The
object $M(\sigma )$ is called the $\sigma$-suspension of $M$. It
is known that $R-gr$ is a Grothendieck category with a family of
projective generators $\{ R(\sigma )|\sigma \in G\}$ (see
\cite{nvo} for details).

Let $I$ be a graded ideal of $R$, and denote ${\cal P}_I=\{ M\in
R-gr|IM=0\}$. Then ${\cal P}_I$ is a closed subcategory of $R-gr$.
Moreover, ${\cal P}_I$ is rigid, i.e. if $M\in {\cal P}_I$, then
$M(\sigma )\in {\cal P}_I$ for any $\sigma \in G$. The smallest
localizing subcategory of $R-gr$ containing ${\cal P}_I$ is
$${\cal C}_I=\{ M\in R-gr|\; {\rm For}\; {\rm any}\; M'<_{R-gr}M, M/M'
\; {\rm contains}\; {\rm some}\; S\in {\cal P}_I, S\neq 0\}$$
Clearly ${\cal C}_I$ is also a rigid subcategory. As in (i) we
obtain that if $M\in R-gr$, then $M$ is ${\cal C}_I$-closed if and
only if the natural morphism
$$\al :M\ra HOM_R(I,M),\;\; \al (m)(a)=am,\; a\in I, m\in M$$
is an isomorphism of graded modules. Recall that
$HOM_R(I,M)=\oplus _{\sigma \in G}HOM_R(I,M)_{\sigma}$, where $
HOM_R(I,M)_{\sigma}$ is the set of all linear maps of degree
$\sigma$ (see \cite{nvo} for details). It is clear from this that
if $M$ is ${\cal C}_I$-closed, then $M(\sigma)$ is ${\cal
C}_I$-closed for any $\sigma \in G$.

(iii)  Let $R$ be a ring with local units, i.e. for any finite
subset $X$ of $R$, there exists an idempotent element $e\in R$
such that $ex=xe=x$ for any $x\in X$ (or equivalently $X\subseteq
eRe$); see \cite{am} for details. For such an $R$ we have the
Grothendieck category $R-MOD$, of all unital left $R$-modules. An
$R$-module $M$ is unital if
$RM=M$. \\
For any two-sided ideal $I$ of $R$, we can define a localizing
category ${\cal C}_I$ as in (i). Also an object $M\in R-MOD$ is
${\cal C}_I$-closed if and only if the natural morphism
$$\al :M\ra RHom_R(I,M),\;\; \al (m)(a)=am,\; a\in I, m\in M$$
is an isomorphism.
\end{ex}

The following is a wide generalization of the classical result of
Kato and M\"{u}ller (see \cite{kato0}, \cite{muller}), which is
given for categories of modules.

\begin{te}  \label{general}
Let $\Gamma =(F,G,\eta ,\rho )$ be a right wide Morita context
between the Grothendieck categories $\cal A$ and $\cal B$, such
that the functors $F$ and $G$ have right adjoints $F'$ and $G'$.
Denote by ${\cal C}_{\Gamma}$ (respectively ${\cal D}_{\Gamma}$)
the subcategory of $\cal A$ (respectively $\cal B$) defined by the
natural morphism $\eta$ (respectively $\rho$). Then the quotient
categories ${\cal A}/\overline{{\cal C}_{\Gamma}}$ and ${\cal
B}/\overline{{\cal D}_{\Gamma}}$ are equivalent via the functors
$F',G'$.
\end{te}
{\bf Proof:} Proposition \ref{propright} shows that ${\cal
C}_{\Gamma}$ is a closed subcategory of $\cal A$, and ${\cal
D}_{\Gamma}$ is a closed subcategory of $\cal B$. By Theorem
\ref{caractclosed}, an object $M\in {\cal A}$ is ${\cal
C}_{\Gamma}$-closed if and only if it is $\overline{{\cal
C}_{\Gamma}}$-closed (and similarly for ${\cal D}_{\Gamma}$-closed
objects). By \cite[pages 195 and 213]{stenstrom}, the category of
$\overline{{\cal C}_{\Gamma}}$-closed objects is equivalent to the
quotient category ${\cal A}/\overline{{\cal C}_{\Gamma}}$, and
this ends the proof. \qed

The following result is a direct consequence of Corollary
\ref{corinj}.

\begin{co}  \label{oneepi}
If $\Gamma =(F,G,\eta ,\rho )$ is a right wide Morita context
between the Grothendieck categories $\cal A$ and $\cal B$, such
that the functors $F$ and $G$ have right adjoints $F'$ and $G'$,
and $\eta$ is an epimorphism, then the categories $\cal A$ and
${\cal B}/\overline{{\cal D}_{\Gamma}}$ are equivalent via the
functors $F'$ and the functor induced by $F$.
\end{co}

\begin{re}
Let $\Gamma =(F,G,\eta ,\rho )$ is a right wide Morita context
between the Grothendieck categories $\cal A$ and $\cal B$.
Proposition \ref{propright}(1) shows that ${\cal D}_{\Gamma}=\{
Y\in {\cal B}|\rho (Y)=0\}$. If we take $Y\in {\rm Ker}\; F$, i.e.
$F(Y)=0$, then $(G\circ F)(Y)=0$, so $\rho (Y)=0$, showing that
$Y\in {\cal D}_{\Gamma}$. Thus for any right wide Morita context
we have ${\rm Ker}\; F\subseteq {\cal D}_{\Gamma}$.

Assume now that $\eta$ is an epimorphism, and let $Y\in {\cal
D}_{\Gamma}$. Then $\eta (F(Y))=F(\rho (Y))=0$. Since $\eta$ is in
fact a natural equivalence, we must have $F(Y)=0$, so $Y\in {\rm
Ker}\; F$. Therefore ${\rm Ker}\; F={\cal D}_{\Gamma}$. Thus
Corollary \ref{oneepi} can be reformulated by saying that $F$
induces an equivalence between the categories ${\cal
B}/\overline{{\rm Ker}\; F}$ and ${\cal A}$.
\end{re}

\section{Applications} \label{sectionapplications}

In this section we apply the general equivalence results that we
proved to several particular cases.

\subsection{Morita contexts for rings with identity and the Kato-M\"{u}ller Theorem}
\label{subrings}

Let $R$ and $S$ be two rings with identity. A Morita context
connecting $R$ and $S$ is a datum $(R,S, _RM_S, _SN_R,\phi ,\psi
)$, where $M$ is an $R$-$S$-bimodule, $N$ is a $S$-$R$-bimodule,
$\phi :M\ot _SN\ra R$ is a morphism of $R$-$R$-bimodules and $\psi
:N\ot _RM\ra S$ is a morphism of $S$-$S$-bimodules such that
\begin{equation} \label{morita1}
\phi (m\ot n)m'=m\psi (n\ot m')
\end{equation}
\begin{equation}  \label{morita2}
 \psi (n\ot m)n'=n\phi (m\ot
n')
\end{equation}
for any $m,m'\in M$, $n,n'\in N$. \\
To such a context we associate two trace ideals: $I={\rm
Im}\;\phi$, which is an ideal of $R$, and $J={\rm Im}\; \psi$,
which is an ideal of $S$. Consider the categories ${\cal A}=R-mod$
and ${\cal B}=S-mod$, and define the functors $F:{\cal B}\ra {\cal
A}$, $F(Y)=M\ot _SY$, and $G:{\cal A}\ra {\cal B}$, $G(X)=N\ot
_RX$. We have a natural morphism $\eta :F\circ G\ra 1_{R-mod}$,
defined by
$$\eta (X):M\ot _SN\ot _RX\ra X,\;\; \eta (X)(m\ot n\ot x)=\phi
(m\ot n)x$$ for any $X\in R-mod$.

We also have a natural morphism $\rho :G\circ F\ra 1_{S-mod}$
defined by
$$\rho (Y):N\ot _RM\ot _SY\ra Y,\;\; \rho (Y)(n\ot m\ot y)=\psi
(n\ot m)y$$ It is straightforward to check that $\Gamma =(F,G,\eta
,\rho )$ is a right wide Morita context. With the notation of
Section \ref{relinj} we have that $I_M={\rm Im}\; \eta (M)=({\rm
Im}\; \phi )M=IM$ for any $M\in R-mod$. By using the definition,
the closed subcategory ${\cal C}_{\Gamma}$ associated to $\Gamma$,
is exactly the subcategory ${\cal P}_I=\{ M\in R-mod|IM=0\}$, so
the smallest localizing subcategory which contains ${\cal
C}_{\Gamma}$ is ${\cal C}_I$ (see Example \ref{exemple} (i)).
Since the functors $F$ and $G$ have right adjoints
$$F':R-mod\ra S-mod,\;\; F'(X)=Hom_R(M,X)$$
$$G':S-mod\ra R-mod,\;\; G'(Y)=Hom_S(N,Y)$$
we see then by Theorem \ref{general} that the quotient categories
$R-mod/{\cal C}_I$ and $S-mod/{\cal C}_J$ are equivalent via the
functors $F'$ and $G'$. This is exactly the Kato-M\"{u}ller
Theorem.

In the case where one of the two maps in the Morita context is
surjective, we obtain the following result (see \cite[Proposition
3.8]{crw}).

\begin{co}  \label{oneepirings}
If $(R,S, _RM_S, _SN_R,\phi ,\psi )$ is a Morita context such that
$\phi$ is surjective, then $R-mod$ is equivalent to a quotient
category of $S-mod$. More precisely, the categories $R-mod$ and
$S-mod/{\cal C}_J$ are equivalent via the functor induced by $F$.
\end{co}
{\bf Proof:} Since $\phi$ is surjective we have that $\eta$ is an
epimorphism. Now we apply Corollary \ref{oneepi}. \qed

The next result shows that in a special case any right wide Morita
context between two categories of modules arises from a Morita
context as we explained above in this subsection.

\begin{pr}
Let $\Delta =(P,Q,\al ,\beta )$ be a right wide Morita context
between the categories ${\cal A}=R-mod$ and ${\cal B}=S-mod$,
where $R$ and $S$ are rings with identity, and assume that the
functors $P$ and $Q$ commute with direct sums. Then there exists a
Morita context $(R,S,M,N,\phi ,\psi )$ connecting $R$ and $S$ such
that $\Delta$ is isomorphic to the right wide Morita context
$\Gamma$ defined by this Morita context.
\end{pr}
{\bf Proof:} Since $P$ is right exact and commutes with direct
sums, there exists an $R,S$-bimodule $M$ such that $P\simeq F=M\ot
_S-$. Similarly there is an $S,R$-bimodule $N$ such that $Q\simeq
G=N\ot _R-$. Let $u:P\ra F$ and $v:Q\ra G$ natural equivalences.
Then $Pv\circ uG:P\circ Q\ra F\circ G$ and $Gu\circ vP:Q\circ P\ra
G\circ F$ are natural equivalences, so there exist natural
transformations $\eta :F\circ G\ra 1_{\cal A}$ and $\rho :G\circ
F\ra 1_{\cal B}$ such that
\begin{equation} \label{eq1isomoritacontext}
\eta \circ Fv\circ uQ=\al
\end{equation}
and
\begin{equation} \label{eq2isomoritacontext}
\rho \circ Gu\circ vP=\beta
\end{equation}
The natural transformation $\eta :F\circ G\ra 1_{\cal A}$ must be
of the form $\eta (X)(m\ot n\ot x)=\phi (m\ot n)x$ for any $X\in
{\cal A}$, $m\in M$, $n\in N$ and $x\in X$, and similarly $\rho
(Y)(n\ot m\ot y)=\psi (n\ot m)y$ for any $Y\in {\cal B}$, $n\in
N$, $m\in M$, $y\in Y$. Moreover, the conditions $P\beta =\al P$
and $Q\al =\beta Q$ imply that $(R,S,M,N,\phi ,\psi )$ is a Morita
context connecting $R$ and $S$. Then $\Gamma =(F,G,\eta ,\rho )$
is the right wide Morita context associated to the Morita context
$(R,S,M,N,\phi ,\psi )$. Moreover, equations
(\ref{eq1isomoritacontext}) and (\ref{eq2isomoritacontext}) show
that $u$ and $v$ give an isomorphism between $\Delta$ and
$\Gamma$. \qed

\subsection{Morita contexts for graded rings}
\label{subgradedrings}

Let $R=\oplus _{\sigma \in G}R_{\sigma}$ and $S=\oplus _{\sigma
\in G}S_{\sigma}$ be two $G$-graded rings, where $G$ is a group. A
graded Morita context is a datum $(R,S, _RM_S, _SN_R,\phi ,\psi
)$, where $M$ is a graded $R$-$S$-bimodule, $N$ is a graded
$S$-$R$-bimodule, $\phi :M\ot _SN\ra R$ is a morphism of graded
$R$-$R$-bimodules and $\psi :N\ot _RM\ra S$ is a morphism of
graded $S$-$S$-bimodules such that equations (\ref{morita1}) and
(\ref{morita2}) are satisfied (see \cite{nvo}). In this case, the
trace ideals of the context, $I={\rm Im}\; \phi$ and $J={\rm Im}\;
\psi$ are graded two-sided ideals. To this graded Morita context
we associate the right wide Morita context $\Gamma =(F,G,\eta
,\rho )$, where
$$F:S-gr\ra R-gr,\;\; F(Y)=M\ot _SY$$
$$G:R-gr\ra S-gr,\;\; G(X)=N\ot _RX$$
and the morphisms $\eta$ and $\rho$ are given by the same formulas
as in Subsection \ref{subrings}. The localizing subcategories
associated to $\Gamma$ are ${\cal C}_I\subseteq R-gr$ and ${\cal
C}_J\subseteq S-gr$ (see Example \ref{exemple} (ii)).

Since the right adjoint functor of $F$ is $F':R-gr\ra S-gr$,
$F'(X)=HOM_R(M,X)$, and the right adjoint of $G$ is $G':S-gr\ra
R-gr$, $G'(Y)=HOM_S(N,Y)$, then we obtain by Theorem \ref{general}
that the quotient categories $R-gr/{\cal C}_I$ and $S-gr/{\cal
C}_J$ are equivalent via the functors $F'$ and $G'$.

\subsection{Morita contexts for rings with local units}
\label{sublocalunits}

Let $R$ and $S$ be two rings with local units, and $R-MOD$ and
$S-MOD$ the associated categories of unital modules. A Morita
context for $R$ and $S$ is a datum $(R,S, _RM_S, _SN_R,\phi ,\psi
)$ as in Subsection \ref{subrings}, with the condition that $M$
and $N$ are unital modules to the left and to the right. The
tensor product is defined exactly as for rings with identity (see
\cite{am} for details). As for rings with identity we obtain by
Theorem \ref{general} an equivalence between the quotient
categories $R-MOD/{\cal C}_I$ and $S-MOD/{\cal C}_J$ via the
functors $SHom_R(M,-)$ and $RHom_S(N,-)$, where $I={\rm Im}\;
\phi$ and $J={\rm Im}\; \psi$.

Also one can obtain a version of Corollary \ref{oneepirings}
exactly in the same way. Note that for instance this explains from
a general point of view \cite[Proposition 3.7]{bdr}.

\subsection{Morita contexts and $I$-projective modules}

Let us consider a Morita context $(R,S, _RM_S, _SN_R,\phi ,\psi )$
connecting the rings with identity $R$ and $S$, and the right wide
Morita context $\Gamma =(F,G,\eta ,\rho )$ associated as in
Subsection \ref{subrings}. Let $F'$ and $G'$ be the right adjoints
of $F$ and $G$ described in Subsection \ref{subrings}, and let
$\eta ':1_{R-mod}\ra G'\circ F'$ be defined as follows
$$\eta '(X):X\ra (G'\circ F')(X)=Hom_S(N,Hom_R(M,X)), \;\; \eta
'(X)(x)(n)(m)=\phi (m\ot n)x$$ for any $X\in R-mod$, $x\in X$,
$m\in M$ and $n\in N$. Similarly one defines $\rho ':1_{S-mod}\ra
F'\circ G'$. In this way we obtain a left wide Morita context
$\Gamma ^{op}=(F',G',\eta ',\rho ')$ between the categories
$R-mod$ and $S-mod$. We call $\Gamma^{op}$ the opposite of
$\Gamma$. If $I={\rm Im}\; \phi$ and $J={\rm Im}\; \psi$, then
${\cal C}_{\Gamma}={\cal C}_{\Gamma ^{op}}={\cal C}_I$ and ${\cal
D}_{\Gamma}={\cal D}_{\Gamma ^{op}}={\cal C}_J$.

If $P\in R-mod$, then $P$ is ${\cal C}_{\Gamma ^{op}}$-projective
if and only if $P$ is $I$-projective, i.e. for any epimorphism
$u:M\ra M'$ with $I{\rm Ker}\; u=0$, and any morphism $f:P\ra M'$,
there exists $g:P\ra M$ such that $u\circ g=f$. \\
On the other hand, $P$ is ${\cal C}_I$-cotorsion free if and only
if $Hom(P,M)=0$ whenever $IM=0$, and it is also equivalent to the
fact that $P=IP$. If we denote by
$${\cal C}_{I,proj}=\{ M\in R-mod|\; M\;{\rm is}\; I-{\rm
projective}\; {\rm and}\; IM=M\}$$
$${\cal C}_{J,proj}=\{ N\in S-mod|\; N\;{\rm is}\; J-{\rm
projective}\; {\rm and}\; JN=N\}$$ then we see by Theorem
\ref{teproj} that the categories ${\cal C}_{I,proj}$ and ${\cal
C}_{J,proj}$ are equivalent via the functors $F$ and $G$.

Similar results can be obtained for the graded case and the local
units case by using the opposite of the right wide Morita contexts
defined in Subsections \ref{subgradedrings} and
\ref{sublocalunits}.

\begin{re}
If $I$ is a two-sided ideal of a ring $R$, the concept of an
$I$-flat module is defined in \cite{kato} as follows. The left
$R$-module $M$ is called $I$-flat if for any exact sequence
$$0\lra N'\stackrel{u}{\lra}N\lra {\rm Coker}\; u\lra 0$$
of right $R$-modules such that $({\rm Coker} \; u)I=0$, we have
that the sequence of abelian groups
$$0\lra N'\ot _RM\stackrel{u\ot 1_M}{\lra}N\ot _RM\lra {\rm Coker}\; u\ot _RM\lra 0$$
is exact. Then one defines a category
$${\cal C}_{I,flat}=\{ M\in R-mod|\; M\;{\rm is}\; I-{\rm
flat}\; {\rm and}\; IM=M\}$$ It is proved in \cite{kato} that if
the $R$-module $M$ is $I$-projective, then it is $I$-flat, and
also that ${\cal C}_{I,proj}={\cal C}_{I,flat}$. An equivalence
result concerning the category ${\cal C}_{I,flat}$, which is a
particular case of our Theorem \ref{teproj}, is proved in
\cite{kato}.
\end{re}

\subsection{Applications to Doi-Hopf modules}

We first recall some facts about coactions of Hopf algebras on
algebras. Let $H$ be a Hopf algebra over a field $k$. Let $A$ be a
right $H$-comodule algebra. This means that $A$ is an algebra, a
right $H$-comodule with $H$-coaction given by $a\mapsto \sum
a_0\ot a_1$ for any $a\in A$, and the comodule structure map from
$A$ to $A\ot H$ is an algebra morphism. The subspace of
coinvariants with respect to this coaction is $A^{co H}=\{ a\in
A|\sum a_0\ot a_1=a\ot 1\}$, and it is a subalgebra of $A$.

We say that $M$ is an $(A,H)$- Doi-Hopf module (or simply a
Doi-Hopf module) if $M$ is a left $A$-module and a right
$H$-comodule (with $m\mapsto\sum m_0\ot m_1$), such that
$$\sum(am)_0\ot (am)_1=\sum a_0m_0\ot a_1m_1$$ for any $a\in A$ and
$m\in M$. We denote by ${_A{\cal M}}^H$ the category whose objects
are the Doi-Hopf modules, and in which the morphisms are the maps
which are $A$-linear and $H$-colinear.

Assume that moreover the Hopf algebra $H$ is co-Frobenius,  i.e.
there exists a non-zero left integral $t\in H^*$. In this case the
rational part $H^{*rat}$ of the dual of $H$ is a subring without
identity of the algebra $H^*$, but $H^{*rat}$ has local units. We
can form the smash product $A\# H^{*rat}$, which is $A\ot
H^{*rat}$ as a vector space (and the element $a\ot h^*$ is denoted
by $a\# h^*$), and has the multiplication given by
$$(a\# h^*)(b\# g^*)=\sum ab_0\# (h^*\leftharpoonup b_1) g^*$$
where $\leftharpoonup$ is the usual right action of $H$ on $H^*$.
It is known that the category ${_A{\cal M}}^H$ is isomorphic to
the category $A\# H^{*rat}-MOD$ of left unital $A\#
H^{*rat}$-modules. The $A\# H^{*rat}$-module structure of a
Doi-Hopf module $M$ is given by $(a\# h^*)\cdot m=\sum
h^*(m_1)am_0$. We identify the categories ${_A{\cal M}}^H$ and
$A\# H^{*rat}-MOD$, i.e. we freely regard a Doi-Hopf module as a
unital module over the smash product, and also the other way
around.

A Morita context $(A\#H^{*rat}, A^{coH}, P, Q, [-,-],(-,-) )$
connecting the smash product and the subalgebra of coinvariants
was constructed in \cite{bdr} (see also \cite[Section 6.3]{dnr}).
We describe briefly this context. We note that in the case where
$H$ is finite dimensional, this context is precisely the one of
\cite{cfm}.

The first bimodule is $P=_{A\#H^{*rat}}A_{A^{coH}}$ with the left
module structure coming from the fact that $A$ itself is a
Doi-Hopf module, and the right module structure obtained by
restriction of scalars. The second bimodule is
$Q=_{A^{coH}}A_{A\#H^{*rat}}$ where the left module structure is
the restriction of scalars, and the right $A\#H^{*rat}$-module
structure is defined by $b \leftarrow (a\#h^*) = \sum b_0 a_0
h^*(S^{-1}(b_1 a_1) g)$, where $S$ is the antipode of $H$ (which
is known to be bijective since $H$ is co-Frobenius), and $g$ is
the distinguished group-like element of $H$ (i.e. $g$ is that
group-like element for which the left integrals on $H$ are exactly
the right $g$-integrals on $H$, see \cite[Section 5.5]{dnr} for
details).

The bimodule maps $[-,-]$ and $(-,-)$ are defined by
$$[-,-]: P \otimes Q = A \otimes  _{A^{coH}} A \ra A \#H^{*rat},$$
$$[-,-](a\otimes b) =[a,b]=\sum ab_0 \#t \leftharpoonup b_1,$$
and
$$(-,-): Q \otimes P = A \otimes  _{A \#H^{*rat}}A   \ra
A^{coH},$$
$$(-,-)(a\otimes b)=(a,b)=t\cdot (ab) =\sum t(a_1b_1)a_0 b_0.$$

The associated trace ideals are $I={\rm Im}\; [-,-]$, an ideal of
$A\# H^{*rat}$, and $J={\rm Im}\; (-,-)=t\cdot A$, an ideal of
$A^{coH}$. The map $[-,-]$ is surjective if and only if the
extension $A/A^{coH}$ is $H$-Galois (see \cite[Section 3]{bdr}),
and the map $(-,-)$ is surjective if and only if there exists a
total integral for $A$, i.e. an $H$-comodule map from $H$ to $A$
that maps 1 to 1 (see \cite[Proposition 3.6]{bdr}). We apply the
results of Subsection \ref{sublocalunits} to this particular
Morita context.

We first note the following.

\begin{lm}
Let $M\in A\#H^{*rat}-MOD$. Then $I\cdot M=0$ if and only if
$t\cdot M=0$. In particular, in the case where $A/A^{coH}$ is
$H$-Galois, $t\cdot M=0$ implies that $M=0$.
\end{lm}
{\bf Proof:} Let $a,b\in A$ and $m\in M$. We have \bea \sum
(ab_0\# t\leftharpoonup b_1)\cdot m&=&\sum (t\leftharpoonup
b_1)(m_1)ab_0m_0\\
&=&\sum t((bm)_1)a(bm)_0\\
&=&a(t\cdot (bm))\eea showing that $I\cdot M=0$ if and only if
$t\cdot M=0$. \qed

Thus the closed category ${\cal P}_I$ (we keep the notation of
Subsection \ref{sublocalunits}) is given by \bea {\cal P}_I&=&\{
M\in A\#H^{*rat}-MOD|I\cdot M=0\}\\
&=&\{ M\in A\#H^{*rat}-MOD|t\cdot M=0\}\eea The smallest
localizing subcategory containing ${\cal P}_I$ is \bea{\cal
C}_I&=&\{ M\in A\#H^{*rat}-MOD|\; {\rm For}\; {\rm any}\; M'<M,
M/M' \; {\rm contains}\; {\rm some}\; S\neq 0\; {\rm with}\;
I\cdot S=0\}\\
&=&\{ M\in A\#H^{*rat}-MOD|\; {\rm For}\; {\rm any}\; M'<M, M/M'
\; {\rm contains}\; {\rm some}\; S\neq 0\; {\rm with}\; t\cdot
S=0\}\eea

Also we define ${\cal P}_J=\{ N\in A^{coH}-mod|JN=0\}$, and the
smallest localizing category ${\cal C}_J$ containing ${\cal P}_J$.
Now we have the following result, which is a particular case of
the results described in Subsection \ref{sublocalunits}.

\begin{pr}  \label{echidoi}
For a co-Frobenius Hopf algebra $H$ and a right Hopf comodule
algebra $A$, there is an equivalence between the quotient
categories $_A{\cal M}^H/{\cal C}_I$ and $A^{coH}-mod/{\cal C}_J$
induced by the above Morita context.
\end{pr}

By the remarks above we have that ${\cal C}_I=0$ if and only
$A/A^{coH}$ is $H$-Galois, and ${\cal C}_J=0$ if and only if there
exists a total integral for $A$. Using these facts, we obtain the
following two particular cases of Proposition \ref{echidoi}.

\begin{co}  \label{cogalois}
If $A$ is a right $H$-comodule algebra such that $A/A^{coH}$ is
$H$-Galois, then the category of Doi-Hopf modules $_A{\cal M}^H$
is equivalent to a quotient category of $A^{coH}-mod$.
\end{co}

\begin{co}  \label{cototal} (\cite[Corollary 3.8]{bdr})
If $A$ is a right $H$-comodule algebra such that there exists a
total integral, then the category $A^{coH}$ is equivalent to a
quotient category of  $_A{\cal M}^H$.
\end{co}

We can explain the equivalence from Corollary \ref{cogalois} in a
more precise way. Assume that $A/A^{coH}$ is $H$-Galois, so ${\cal
C}_I=0$. Let
$$F':{_A{\cal M}}^H\ra A^{coH}-mod, \;\; F'(X)=Hom_{A
\#H^{*rat}}(A,X)$$ be the right adjoint of the functor
$$F=A\ot ^{A^{coH}}-:A^{coH}-mod\ra {_A{\cal M}}^H$$
It is easy to see that for any Doi-Hopf module $X$ we have a
natural isomorphism $$Hom_{A \#H^{*rat}}(A,X)=Hom_{_A{\cal
M}^H}(A,X)\simeq X^{coH}$$ Therefore $F'\simeq (-)^{coH}$, the
functor that takes the coinvariants of a Doi-Hopf module.

Let $\cal D$ be the subcategory of $A^{coH}-mod$ consisting of all
${\cal C}_J$-closed objects. Then by Theorem \ref{teinj}, $F'$ is
an equivalence between the categories $_A{\cal M}^H$ and $\cal D$
(the inverse of this equivalence is $G'$, the right adjoint of the
functor $G=A\ot _{A \#H^{*rat}}-$). If we restrict $F$ to $\cal
D$, we see that it is still a left adjoint of the equivalence
functor $F'$ (regarded from $_A{\cal M}^H$ to $\cal D$). We
conclude that when restricted to $\cal D$, the functor $F$ is
itself an equivalence, so $F\circ F'\simeq 1_{_A{\cal M}^H}$. This
means that for any Doi-Hopf module $M$, the natural morphism
$$\phi _M:A\ot _{A^{coH}}M^{coH}\ra M,\;\; \phi _M(a\ot m)=am$$
is an isomorphism of Doi-Hopf modules. This is exactly the Weak
Structure Theorem for Galois extensions. This result is not new.
It was obtained in \cite[Theorem 3.1]{bdr}, and it also follows
from \cite[2.11]{dt} in presence of the fact that for an
$H$-Galois extension $A/A^{coH}$ with $H$-co-Frobenius, the left
$A^{coH}$-module $A$ is flat (see \cite[Corollary 3.5]{bdr}).
However our approach gives a more categorial idea about how the
Weak Structure Theorem arises.

\section{Left wide Morita contexts and an equivalence theorem}
\label{leftwideequivalence}

Let $\Gamma =(F,G,\eta ,\rho )$ be a left wide Morita context
between the Grothendieck categories $\cal A$ and $\cal B$. As in
Section \ref{relproj} we denote by
$${\cal C}_{\Gamma}=\{ M\in {\cal A}|\eta (M)=0\}$$
the closed subcategory associated to $\Gamma$. Let
$$\overline{{\cal C}_{\Gamma}}=\{ X\in {\cal A}|\; {\rm For}\; {\rm any}\; X'<X,
X/X'\; {\rm contains}\; {\rm a}\; {\rm non-zero}\; {\rm object}\;
{\rm of}\; {\cal C}_{\Gamma}\}$$ be the smallest localizing
subcategory that contains ${\cal C}_{\Gamma}$. Similarly we define
the subcategory ${\cal D}_{\Gamma}$ of $\cal B$, and the smallest
localizing subcategory $\overline{{\cal D}_{\Gamma}}$ containing
${\cal D}_{\Gamma}$. We see by the proof of Theorem
\ref{caractclosed} that an object $M\in {\cal A}$ is ${\cal
C}_{\Gamma}$-torsion free if and only if $M$ is $\overline{{\cal
C}_{\Gamma}}$-torsion free.

\begin{lm} \label{keretatf}
If $M\in {\cal A}$ is ${\cal C}_{\Gamma}$-torsion free, then $\eta
(M):M\ra (F\circ G)(M)$ is a monomorphism, and $G(M)$ is ${\cal
D}_{\Gamma}$-torsion free.
\end{lm}
{\bf Proof:} By Proposition \ref{propleft} we have ${\rm Ker} \;
\eta (M)\in {\cal C}_{\Gamma}$. Since $M$ is ${\cal
C}_{\Gamma}$-torsion free, we get that ${\rm Ker}\; \eta (M)=0$,
so $\eta (M)$ is a monomorphism.

Assume that $G(M)$ were not ${\cal D}_{\Gamma}$-torsion free. Then
there exists $Y\in {\cal D}_{\Gamma}$, $Y\neq 0$ such that $Y\leq
G(M)$. Let $i:Y\ra G(M)$ be the inclusion morphism. Since $Y\in
{\cal D}_{\Gamma}$, we have $\rho (Y)=0$. Then  $\rho (G(M))\circ
i=(G\circ F)(i)\circ \rho (Y)=0$. Since $\rho G=G\eta$ we obtain
that $G\eta (M)\circ i=0$. But $\eta (M)$ is a monomorphism and
$G$ is left exact, so $G\eta (M)$ is also a monomorphism. This
shows that $i=0$, so $Y=0$, a contradiction. We conclude that
$G(M)$ must be ${\cal D}_{\Gamma}$-torsion free. \qed

\begin{lm}  \label{gcgamma}
(i) If $M\in {\cal C}_{\Gamma}$, then $G(M)\in {\cal
D}_{\Gamma}$.\\
(ii) If the functors $F$ and $G$ commute with direct limits, then
for any $M\in \overline{{\cal C}_{\Gamma}}$ we have that $G(M)\in
\overline{{\cal D}_{\Gamma}}$.
\end{lm}
{\bf Proof:} (i) Since $M\in {\cal C}_{\Gamma}$ we have $\eta
(M)=0$. Then $\rho (G(M))=G(\eta (M))=0$, so $G(M)\in M\in {\cal
D}_{\Gamma}$.\\
(ii) Let $t$ be the preradical associated to the closed
subcategory ${\cal C}_{\Gamma}$. We define by transfinite
recurrence the objects $M_{\al}$ for any ordinal $\al$. We first
set $M_1=t(M)\in {\cal C}_{\Gamma}$. If $\al$ is an ordinal such
that $M_{\al}$ is defined, we define $M_{\al +1}$ such that
$M_{\al +1}/M_{\al}=t(M/M_{\al})$. Finally, if $\al$ is a limit
ordinal, we put $M_{\al}=\cup _{\beta <\al}M_{\beta}$.

Since $M\in \overline{{\cal C}_{\Gamma}}$, there exists an ordinal
$\al _0$ such that $M=M_{\al _0}$. We prove by transfinite
induction that $G(M_{\al})\in \overline{{\cal D}_{\Gamma}}$ for
any ordinal $\al \leq \al _0$. Indeed, for $\al =1$ it follows by
the assertion (i). If $\al$ is a limit ordinal, we have that
$$G(M_{\al})=G(\cup _{\beta <\al}M_{\beta})=\lim_{\stackrel{\longrightarrow}{\beta
<\al}}\; G(M_{\beta})\in \overline{{\cal D}_{\Gamma}}$$ since $G$
commutes with direct limits and $\overline{{\cal D}_{\Gamma}}$ is
a localizing subcategory. \\
If $\al $ is an arbitrary ordinal, we have the exact sequence
$$0\lra M_{\al}\lra M_{\al +1}\stackrel{\pi_{\al +1}}{\lra}M_{\al
+1}/M_{\al} {\lra} 0$$ where $M_{\al}\in \overline{{\cal
C}_{\Gamma}}$ and $M_{\al +1}/M_{\al}\in {\cal C}_{\Gamma}\subset
\overline{{\cal C}_{\Gamma}}$. Since $G$ is left exact, we have
the exact sequence
$$0\lra G(M_{\al})\lra G(M_{\al +1})\stackrel{G(\pi_{\al +1})}{\lra}G(M_{\al
+1}/M_{\al})$$ Since $G(M_{\al +1}/M_{\al})\in {\cal D}_{\Gamma}$,
we have that ${\rm Im}\; G(\pi _{\al +1})\in {\cal D}_{\Gamma}$.
This shows that $G(M_{a\ +1})\in \overline{{\cal D}_{\Gamma}}$,
since $\overline{{\cal D}_{\Gamma}}$ is closed under
extensions.\qed

\begin{pr}  \label{etaiso}
If $M$ is $\overline{{\cal C}_{\Gamma}}$-closed, then $\eta
(M):M\ra (F\circ G)(M)$ is an isomorphism. If $N\in {\cal B}$ is
$\overline{{\cal D}_{\Gamma}}$-closed, then $\rho (N):N\ra (G\circ
F)(N)$ is an isomorphism.
\end{pr}
{\bf Proof:} Since $M$ is $\overline{{\cal C}_{\Gamma}}$-torsion
free, it is also ${\cal C}_{\Gamma}$-torsion free, so by Lemma
\ref{keretatf} we have that $\eta (M)$ is a monomorphism and
$(F\circ G)(M)$ is ${\cal C}_{\Gamma}$-torsion free. On the other
hand ${\rm Coker}\; \eta (M)\in {\cal C}_{\Gamma}$ by Proposition
\ref{cokercgamma}, so $\eta (M)$ is an essential monomorphism.
Since $M$ is $\overline{{\cal C}_{\Gamma}}$-injective, there
exists a morphism $g:(F\circ G)(M)\ra M$ such that $g\circ \eta
(M)=1_M$. This shows that $(F\circ G)(M)\simeq {\rm Im}\; \eta
(M)\oplus {\rm Coker}\; \eta (M)$. But $(F\circ G)(M)$ is
$\overline{{\cal C}_{\Gamma}}$-torsion free, so ${\rm Coker} \;
\eta (M)$ must be 0, and then $\eta (M)$ is an isomorphism.
Similarly we see that if $N\in {\cal B}$ is $\overline{{\cal
D}_{\Gamma}}$-closed, then $\rho (N):N\ra (G\circ F)(N)$ is an
isomorphism. \qed

\begin{pr}  \label{imclosed}
If $M\in {\cal A}$ is $\overline{{\cal C}_{\Gamma}}$-closed, then
$G(M)$ is $\overline{{\cal D}_{\Gamma}}$-closed.  Similarly, if
$N\in {\cal B}$ is $\overline{{\cal D}_{\Gamma}}$-closed, then
$F(M)$ is $\overline{{\cal C}_{\Gamma}}$-closed.
\end{pr}
{\bf Proof:} Assume that $M\in {\cal A}$ is $\overline{{\cal
C}_{\Gamma}}$-closed. By Lemma \ref{keretatf} we have that $G(M)$
is ${\cal D}_{\Gamma}$-torsion free, so it is $\overline{{\cal
D}_{\Gamma}}$-torsion free, too. Let $N\in {\cal B}$ be the
closure of $G(M)$ with respect to the localizing subcategory
$\overline{{\cal D}_{\Gamma}}$. This means that if $T':{\cal B}\ra
{\cal B}/\overline{{\cal D}_{\Gamma}}$ is the natural functor
associated to the quotient category ${\cal B}/\overline{{\cal
D}_{\Gamma}}$, and $S'$ is a right adjoint of $T'$, then
$N=(S'\circ T')(G(M))$. We have the exact sequence
$$0\lra G(M)\stackrel{i}{\lra}N\lra {\rm Coker}\; i\lra 0
$$ with ${\rm Coker} \; i\in \overline{{\cal D}_{\Gamma}}$. Since $F$ is
left exact, we obtain the exact sequence
$$0\lra (F\circ G)(M)\stackrel{F(i)}{\lra}F(N)\stackrel{\theta}{\lra}F({\rm Coker}\; i)
$$ with $F({\rm Coker} \; i)\in \overline{{\cal C}_{\Gamma}}$, by Lemma
\ref{gcgamma}(ii). We have that ${\rm Im}\; \theta \in
\overline{{\cal D}_{\Gamma}}$. By Proposition \ref{etaiso} we have
$(F\circ G)(M)\simeq M$, so $(F\circ G)(M)$ is $\overline{{\cal
D}_{\Gamma}}$-closed. Now the diagram

\begin{picture}(100,110)(-30,-10)
\put(48,47){$0$} \put(60,51){\vector(1,0){38}}\put(82,53){$i$}
 \put(100,47){$(F\circ G)(M)$}
\put(160,51){\vector(1,0){50}}\put(171,53){$F(i)$}
\put(215,47){$F(N)$}\put(250,51){\vector(1,0){30}}
\put(261,53){$j$}\put(285,47){${\rm Im}\; \theta$}
\put(128,43){\vector(0,-1){40}} \put(100,-10){$(F\circ G)(M)$}
\put(132,20){$1_{(F\circ G)(M)}$} \put(315,51){\vector(1,0){20}}
\put(338,47){$0$}
\end{picture}

shows that there exists a morphism $h:F(N)\ra (F\circ G)(M)$ such
that $h\circ F(i)=1_{(F\circ G)(M)}$. Hence $F(N)\simeq {\rm Im}\;
F(i)\oplus {\rm Im}\; \theta$. Since $F(N)$ is $\overline{{\cal
C}_{\Gamma}}$-closed, then it is $\overline{{\cal
C}_{\Gamma}}$-torsion free, so ${\rm Im}\; \theta =0$. This shows
that $F(i)$ is an isomorphism. We have the morphisms

\begin{picture}(100,120)(-40,20)
\put(92,120){$(G\circ F\circ G)(M)$}
\put(178,123){\vector(1,0){75}}\put(190,126){$(G\circ F)(i)$}
 \put(257,120){$(G\circ F)(N)$}
\put(257,50){$N$} \put(261,60){\vector(0,1){53}}\put(264,90){$\rho
(N)$} \put(17,123){\vector(1,0){72}}\put(30,126){$G(\eta (M))$}
\put(3,120){$M$}
\end{picture}

where $\eta (M)$ is an isomorphism and $\rho (N)$ is also an
isomorphism by Proposition \ref{etaiso}. We conclude that
$G(M)\simeq N$, so $G(M)$ is $\overline{{\cal
D}_{\Gamma}}$-closed. \qed

\begin{te}  \label{ecleft}
Let $\Gamma =(F,G,\eta ,\rho )$ be a left wide Morita context
between the Grothendieck categories $\cal A$ and $\cal B$ such
that the functors $F$ and $G$ commute with direct limits. Then the
quotient categories ${\cal A}/\overline{{\cal C}_{\Gamma}}$ and
${\cal B}/\overline{{\cal D}_{\Gamma}}$ are equivalent via the
restriction of the functors $F$ and $G$.
\end{te}
{\bf Proof:} Since ${\cal A}/\overline{{\cal C}_{\Gamma}}$ is the
subcategory of all $\overline{{\cal C}_{\Gamma}}$-closed objects
of $\cal A$, and ${\cal B}/\overline{{\cal D}_{\Gamma}}$ is the
subcategory of all $\overline{{\cal D}_{\Gamma}}$-closed objects
of $\cal B$, the result follows from Propositions \ref{etaiso} and
\ref{imclosed}. \qed

\section{Applications to Morita-Takeuchi contexts}

In this section we apply the results of Section
\ref{leftwideequivalence} to left wide Morita contexts arising
from Morita-Takeuchi contexts.

\subsection{Morita-Takeuchi contexts and a theorem of Berbec}

Let $C$ and $D$ be two coalgebras over a field. A Morita-Takeuchi
context connecting $C$ and $D$ is a datum $(C,D,M,N,\phi ,\psi )$,
where $M$ is a $C,D$-bicomodule, $N$ is a $D,C$-bicomodule, $\phi
:C\ra M\Cut _D N$ is a morphism of $C,C$-bicomodules, and $\psi
:D\ra N\Cut _CM$ is a morphism of $D,D$-bicomodules such that
$(1_M\Cut \psi )\circ \gamma _{M,D}=(\phi \Cut 1_M)\circ \gamma
_{C,M}$ and $(1_N\Cut \phi )\circ \gamma _{N,C}=(\psi \Cut
1_N)\circ \gamma _{D,N}$, where $\Cut$ is the cotensor product,
$\gamma_{M,D}:M\ra M\Cut _DD$, $\gamma_{C,M}:M\ra C\Cut _CM$,
$\gamma _{N,C}:N\ra N\Cut _CC$ and $\gamma_{D,N}:N\ra D\Cut _DN$
are the natural isomorphisms.

To such a Morita-Takeuchi context we associate a left wide Morita
context $\Gamma =(F,G,\eta ,\rho )$ between the categories of
right comodules ${\cal M}^D$ and ${\cal M}^C$ (see \cite{cg}),
where the functors $F$ and $G$ are defined by
$$F:{\cal M}^C\ra {\cal M}^D,\;\; F(X)=X\Cut _CM$$
$$G:{\cal M}^D\ra {\cal M}^C,\;\; G(Y)=Y\Cut _DN$$
and the natural morphisms $\eta :1_{{\cal M}^D}\ra F\circ G$ and
$\rho :1_{{\cal M}^C}\ra G\circ F$ are defined by
$$\eta (Y)=(1_Y\Cut _D\psi )\circ \gamma _{Y,D} \;\; {\rm for}\; {\rm any}\; Y\in {\cal M}^D$$
$$\rho (X)=(1_X\Cut _C\phi )\circ \gamma _{X,C} \;\; {\rm for}\; {\rm any}\; X\in {\cal M}^C$$
The functors $F$ and $G$ are left exact and commute with direct
limits. Denote $A={\rm Ker}\; \phi$, which is a subcoalgebra of
$C$, and $B={\rm Ker}\; \psi$, which is a subcoalgebra of $D$. Let
${\cal C}_{\Gamma}$ be the closed subcategory of ${\cal M}^C$
defined by $\Gamma$ as in Section \ref{relinj}. We have that \bea
{\cal
C}_{\Gamma}&=&\{ M\in {\cal M}^C|M\Cut _CA\simeq M\}\\
&=&\{ M\in {\cal M}^C|\rho _M(M)\subseteq M\ot A\}\\
&=&\{ M\in {\cal M}^C|A^{\perp}M=0\}\eea where $A^{\perp}$ is the
subspace of $C^*$ consisting of all maps that annihilates $A$ (see
\cite[Sections 1.2 and 2.5]{dnr}). The smallest localizing
subcategory containing ${\cal C}_{\Gamma}$ is
$$\overline{{\cal C}_{\Gamma}}=\{ M\in {\cal
M}^C|A^{\perp}_{\infty}M=0\}$$ where $A_{\infty}=\cup _{n\geq 1}
\wedge ^nA$ (here $\wedge$ is the usual wedge, see \cite{dnr},
\cite{mo}).

Now we can derive in a natural way as a particular case of Theorem
\ref{general} the following result of Berbec (see \cite{berbec}).

\begin{co}  \label{teoremaberbec}
Let $(C,D,M,N,f,g)$ be a Morita-Takeuchi context connecting the
coalgebras $C$ and $D$, and let $\Gamma$ be the associated left
wide Morita context as above. Then the quotient categories ${\cal
M}^C/\overline{{\cal C}_{\Gamma}}$ and ${\cal M}^D/\overline{{\cal
D}_{\Gamma}}$ are equivalent.
\end{co}
{\bf Proof:} The cotensor product functors are left exact and
commute with direct limits, so the result follows directly from
Theorem \ref{ecleft}. \qed

In the particular case where one of the maps of the
Morita-Takeuchi context is injective, we obtain the following
result (see \cite[Proposition 2.2]{dnrv}).

\begin{co}  \label{mtcontextfinj}
Let $(C,D,M,N,f,g)$ be a Morita-Takeuchi context such that $f$ is
injective. Then the category ${\cal M}^C$ is equivalent to a
quotient category of ${\cal M}^D$.
\end{co}

\subsection{Applications to Hopf-Galois coextensions}

Let $H$ be a finite dimensional Hopf algebra over a field $k$, and
let $C$ be a left $H$-comodule coalgebra. This means that $C$ is a
coalgebra (with comultiplication $c\mapsto \sum c_1\ot c_2$) and a
left $H$-comodule (with coaction $c\mapsto \sum c_{(-1)}\ot
c_{(0)}$) such that
$$\sum c_{(-1)}\ot c_{(0)1}\ot c_{(0)2}=\sum c_{1(-1)}c_{2(-1)}\ot
c_{1(0)}\ot c_{2(0)}$$ and $$\sum \varepsilon
(c_{(0)})c_{(-1)}=\varepsilon (c)1_H$$ for any $c\in C$. We can
form the smash coproduct $C\tie H$, which is $C\ot H$ as a
$k$-vector space, with the element $c\ot h$ denoted by $c\tie h$,
and has a coalgebra structure with counit $\varepsilon _C\tie
\varepsilon _H$ and comultiplication given by
$$\Delta (c\tie h)=\sum (c_1\tie c_{2(-1)}h_1)\ot (c_{2(0)}\tie
h_2)$$ We also consider the factor coalgebra
$\overline{C}=C/CH^{*+}$, where $C$ is regarded as a right
$H^*$-module, and $H^{*+}={\rm Ker}\; \varepsilon _{H^*}$. Let
$t\in H^*$ be a left integral on $H$, and $a\in H$ the
distinguished grouplike element. Then we have a Morita-Takeuchi
context $(C\tie H,\overline{C},C,C,f,g)$ as follows (see
\cite[Theorem 1.1]{drz} for details). The left and right
$\overline{C}$-comodule structures on $C$ come via the natural
projection $C\ra \overline{C}$. The left coaction of $C\tie H$ on
$C$ is given by $c\mapsto \sum (c_1\tie c_{2(-1)})\ot c_{2(0)}$,
and the right coaction of $C\tie H$ on $C$ is
$$c\mapsto \sum c_{1(0)}\ot (c_{2(0)}\tie
S^{-1}(c_{1(-1)}c_{2(-1)})a)$$ The maps $f$ and $g$ are defined by
$$f:C\tie H\ra C\Cut _{\overline{C}}C, \; f(c\tie h)=\sum c_1\Cut
t(c_{2(-1)}h)c_{2(0)}$$
$$g:\overline{C}\ra C\Cut _{C\tie H}C,\; g(\overline{c})=\sum
t(c_{1(-1)}c_{2(-1)})c_{1(0)}\Cut c_{2(0)}$$ where $\overline{c}$
denotes the class of $c\in C$ modulo the coideal $CH^{*+}$.

We can apply Corollary \ref{teoremaberbec} to this Morita-Takeuchi
context, and we find that certain quotient categories of ${\cal
M}^{C\tie H}$ and ${\cal M}^{\overline{C}}$ are equivalent. This
may be reformulated if we take into account that the category
${\cal M}^{C\tie H}$ is isomorphic to the category of right
$C,H$-comodules, consisiting of all objects that are right
$C$-comodules and right $H$-comodules, and the two comodule
structures satisfy a compatibility condition (see \cite{cdr}).

If moreover $C/\overline{C}$ is an $H^*$-Galois coextension, which
is equivalent to the map $f$ in the Morita-Takeuchi context being
injective (see \cite[Theorem 1.2]{drz}), then we obtain by
Corollary \ref{mtcontextfinj} that the category ${\cal M}^{C\tie
H}$ is equivalent to a quotient category of ${\cal
M}^{\overline{C}}$.

On the other hand, in the case where $H$ is cosemisimple (or
equivalently $H^*$ is semisimple), we have by \cite[Proposition
3.7 and the comments before it]{cdr} that $\overline{C}\simeq
C^{coH}$, the associated coalgebra of coinvariants. In this case
it is easy to see that the map $g$ is injective. Indeed, if
$g(\overline{c})=0$, then by applying $I\ot \varepsilon _C$, we
get that $c\cdot t=0$. Since $H$ is cosemisimple we can choose $t$
such that $t(1)=1$. Then $(\varepsilon -t)(1)=0$, so $\varepsilon
-t\in H^{*+}$. Hence $c=c\cdot \varepsilon =c \cdot (\varepsilon
-t)\in CH^{*+}$, so $\overline{c}=0$. Thus for cosemisimple $H$ we
obtain by Corollary \ref{mtcontextfinj} that the category ${\cal
M}^{\overline{C}}$ is equivalent to a quotient category of ${\cal
M}^{C\tie H}$.


\begin{thebibliography}{99}

\bibitem{an} T. Albu and C. N\u{a}st\u{a}sescu, Relative
finiteness in module theory, Marcel Dekker, 1984.

\bibitem{af} F.W. Anderson, K.R. Fuller, Rings and Categories of Modules,
 GTM 13,  2-nd edition, Springer-Verlag, 1992.


\bibitem{am}  P.N. \'Anh and L. M\'arki, Morita equivalences for rings
without identity, Tsukuba J. Math. {\bf 11} (1987), 1-16.

\bibitem{bdr} M. Beattie, S. D\u{a}sc\u{a}lescu, and \c{S}.
Raianu, Galois Extensions for Co-Frobenius Hopf Algebras, J.
Algebra {\bf 198}(1997), 164-183.

\bibitem{berbec} I. Berbec, The Morita-Takeuchi theory for
quotient categories, Comm. Algebra {\bf 31} (2003), No.2, 843-858.

\bibitem{bm} R. Blattner, S. Montgomery, A duality theorem for Hopf module algebras,
J. Algebra 95 (1985), 153-172.

\bibitem{cdr} S. Caenepeel, S. D\u{a}sc\u{a}lescu
and  \c{S}. Raianu, Cosemisimple Hopf algebras coacting on
coalgebras, Comm. Algebra {\bf 24} (1996), 1649-1677.

\bibitem{cc} C. Cai and H. Chen, Coactions, smash products, and Hopf modules,
J. Algebra {\bf 167} (1994), 85-99.

\bibitem{cg0} F. Casta\~{n}o Iglesias and J. Gomez Torrecillas,
Wide Morita contexts, Comm. Algebra {\bf 23} (1995), 601-622.

\bibitem{cg} F. Casta\~{n}o Iglesias and J. Gomez Torrecillas,
Wide Morita contexts and equivalences of comodule categories, J.
Pure Appl. Algebra {\bf 131} (1998), 213-225.


\bibitem{cfm} M. Cohen, D. Fischman and S. Montgomery, Hopf
Galois extensions, smash products and Morita equivalence, J.
Algebra {\bf 133} (1990), 351-372.

\bibitem{crw} M. Cohen, \c{S}. Raianu and S. Westreich, Semiinvariants for
Hopf algebra actions, Israel J. Math. 88 (1994), 279-306.

\bibitem{dnrv}
 S. D\u{a}sc\u{a}lescu, C. N\u{a}st\u{a}sescu, \c{S}. Raianu and F. Van Oystaeyen,
Graded coalgebras and Morita-Takeuchi contexts, Tsukuba J. Math.
{\bf 19}(1995), 395-407.

\bibitem{dnr} S. D\u{a}sc\u{a}lescu, C. N\u{a}st\u{a}sescu
and  \c{S}. Raianu, Hopf algebras: an introduction, Pure and
Applied Math. {\bf 235} (2000), Marcel Dekker.


\bibitem{drz} S. D\u{a}sc\u{a}lescu, \c{S}. Raianu, Y. H. Zhang, Finite Hopf-Galois
Coextensions, Crossed Coproducts and Duality, J. Algebra  {\bf
178}(1995), 400-413.

\bibitem{meyeringraham} F. De Meyer and E. Ingraham, Separable
algebras over commutative rings, Springer Verlag, 1971.

\bibitem{dt} Y. Doi, M. Takeuchi, Hopf-Galois extensions of algebras, the
Miyashita-Ulbrich action and Azumaya algebras, J. Algebra {\bf
121} (1989), 488-516.

\bibitem{F} C. Faith, Algebra: Rings, Modules and Categories, I, Springer,
Berlin, 1973.


\bibitem{gabriel} P. Gabriel, Des categories abeliennes, Bull.
Soc. Math. France {\bf 90} (1962), 323-448.

\bibitem{kato0} T. Kato, $U$-distinguished modules, J. Algebra
{\bf 25} (1973), 15-24.

\bibitem{kato} T. Kato and K. Ohtake, Morita context and
equivalence, J. Algebra {\bf 61} (1979), 360-366.

\bibitem{macl} S. Mac Lane, Categories for the working
mathematician, Springer Verlag, 1994.

\bibitem{mo}
S. Montgomery, Hopf algebras and their actions on rings, CBMS Reg.
Conf. Series 82, Providence, R.I., 1993.

\bibitem{muller} B. J. M\"{u}ller, The quotient category of a
Morita context, J. Algebra {\bf 28} (1974), 389-407.

\bibitem{nvo} C. N\u{a}st\u{a}sescu and F. Van Oystaeyen, Graded
ring theory, North Holland, 1982.

\bibitem{nw} W. K. Nicholson and J. F. Watters, Morita context functors, Math.
Proc. Cambridge Phil. Soc. {\bf 103} (1988), 399-408.

\bibitem{stenstrom} B. Stenstr\"{o}m, Rings of quotients, Springer
Verlag, 1975.

\bibitem{ta}
M. Takeuchi,  Morita Theorems for Categories of Comodules,  J.
Fac. of Sci. Univ. Tokyo, {\bf 24} (1977), 629-644.


\bibitem{vdb} M. Van den Bergh, A duality theorem for Hopf algebras, in Methods in Ring
Theory, NATO ASI Series vol. 129, Reidel, Dordrecht, 1984,
517-522.

\bibitem{wisbauer} R. Wisbauer, Foundations of module and ring
theory, Gordon and Breach, 1991.

\end{thebibliography}
\end{document}